\documentclass[11pt,twoside,english]{myamsart}
\advance\oddsidemargin by -1.0cm
\advance\evensidemargin by -1.0cm
\advance\topmargin by -1.0cm 
\usepackage{amssymb}
\usepackage{babel}
\usepackage{amstext}
\usepackage{amscd}   
\usepackage{epsfig}  
\usepackage{rotating}
\let\mathg\mathfrak
\theoremstyle{plain}
\newtheorem{cor}{Corollary}[section]
\newtheorem{lem}{Lemma}[section]
\newtheorem{thm}{Theorem}[section]            
\newtheorem{prop}{Proposition}[section]
\theoremstyle{definition}
\newtheorem{exa}{Example}[section]
\newtheorem{NB}{Remark}[section]

\renewcommand{\labelenumi}{$(\theenumi)$}
\newcommand{\bdm}{\begin{displaymath}}
\newcommand{\edm}{\end{displaymath}}
\newcommand{\be}{\begin{equation}}
\newcommand{\ee}{\end{equation}}
\newcommand{\ba}[1]{\begin{array}{#1}}
\newcommand{\ea}{\end{array}}

\newcommand{\btab}{\begin{tabular}}
\newcommand{\etab}{\end{tabular}}

\newcommand{\C}{\ensuremath{\mathbb{C}}}
\newcommand{\R}{\ensuremath{\mathbb{R}}}
\newcommand{\CP}{\ensuremath{\mathbb{CP}}}

\newcommand{\Z}{\ensuremath{\mathbb{Z}}}

\newcommand{\M}{\ensuremath{\mathcal{M}}}

\newcommand{\T}{\ensuremath{\mathrm{T}}}

\newcommand{\G}{\ensuremath{\mathrm{G}}}
\newcommand{\J}{\ensuremath{\mathrm{J}}}

\newcommand{\vrho}{\ensuremath{\varrho}}

\newcommand{\Scal}{\ensuremath{\mathrm{Scal}}}

\newcommand{\Lin}{\ensuremath{\mathrm{Lin}}}

\newcommand{\kl}[1]{\ensuremath{\scriptstyle{#1}}}

\newcommand{\slin}{\ensuremath{\mathg{sl}}}
\newcommand{\SL}{\ensuremath{\mathrm{SL}}}
\newcommand{\un}{\ensuremath{\mathg{u}}}
\newcommand{\su}{\ensuremath{\mathg{su}}}
\newcommand{\SU}{\ensuremath{\mathrm{SU}}}
\newcommand{\Sp}{\ensuremath{\mathrm{Sp}}}
\newcommand{\U}{\ensuremath{\mathrm{U}}}
\newcommand{\Ni}{\ensuremath{\mathrm{N}}}

\newcommand{\so}{\ensuremath{\mathg{so}}}
\newcommand{\SO}{\ensuremath{\mathrm{SO}}}

\newcommand{\p}{\ensuremath{\mathg{p}}}

\newcommand{\g}{\ensuremath{\mathfrak{g}}}
\newcommand{\h}{\ensuremath{\mathfrak{h}}}
\newcommand{\m}{\ensuremath{\mathfrak{m}}}

\newcommand{\z}{\ensuremath{\mathfrak{z}}}

\begin{document}
\def\haken{\mathbin{\hbox to 6pt{%
                 \vrule height0.4pt width5pt depth0pt
                 \kern-.4pt
                 \vrule height6pt width0.4pt depth0pt\hss}}}
    \let \hook\intprod
\setcounter{equation}{0}
%
%
\thispagestyle{empty}
%
\date{\today}
\title[Almost hermitian $6$-manifolds revisited]
{Almost hermitian $6$-manifolds revisited}
%
%
%
\author{Bogdan Alexandrov}
\author{Thomas Friedrich}
\author{Nils Schoemann}
\address{\hspace{-5mm} 
{\normalfont\ttfamily bogdan@mathematik.hu-berlin.de}\newline
{\normalfont\ttfamily friedric@mathematik.hu-berlin.de}\newline
{\normalfont\ttfamily schoeman@mathematik.hu-berlin.de}\newline
Institut f\"ur Mathematik \newline
Humboldt-Universit\"at zu Berlin\newline
Sitz: WBC Adlershof\newline
D-10099 Berlin, Germany}
%
\thanks{Supported by the SFB 288 "Differential Geometry
and Quantum Physics" and the SPP 1154 ``Globale Differentialgeometrie'' 
of the DFG}
\subjclass[2000]{Primary 53 C 25; Secondary 81 T 30}
\begin{abstract}
A Theorem of Kirichenko states that the torsion 
$3$-form of the characteristic connection of a
nearly K\"ahler manifold is parallel. On the other side, any 
almost hermitian manifold of type $\mathrm{G}_1$ admits a unique
connection with totally skew symmetric torsion.
In dimension six, we generalize Kirichenko's Theorem and 
we describe almost hermitian $\mathrm{G}_1$-manifolds 
with parallel torsion form. In particular,  
among them there are only two types of $\mathcal{W}_3$-manifolds 
with a non-abelian holonomy group, namely
twistor spaces of $4$-dimensional self-dual Einstein manifolds
and the invariant hermitian structure on the Lie group $\mathrm{SL}(2, \C)$. 
Moreover, we classify all naturally reductive hermitian 
$\mathcal{W}_3$-manifolds with  small isotropy group of the
characteristic torsion. 
\end{abstract}
\maketitle
\tableofcontents
\pagestyle{headings}
%
%
%
\section{Introduction}\noindent
\noindent
Fix a subgroup $\mathrm{G} \subset \SO(n)$ of the special
orthogonal group and decompose the Lie algebra $\so(n) = \g \oplus \m$ 
into the Lie algebra $\g$ of $\mathrm{G}$ and its orthogonal complement $\m$. 
The different geometric types of  
$\mathrm{G}$-structures on a Riemannian manifold correspond to the 
irreducible $\mathrm{G}$-components of the representation $\R^n \otimes \m$.
This approach to non integrable geometries is a kind of folklore in
differential geometry, and was exposed in detail in
the article \cite{Fri3}. Indeed, consider a Riemannian manifold 
$(M^n, g)$ and denote 
its Riemannian frame bundle by 
$\mathcal{F}(M^n)$. It is a principal $\SO(n)$-bundle 
over $M^n$. A $\mathrm{G}$-structure is a reduction 
$\mathcal{R} \subset \mathcal{F}(M^{n})$ of the frame bundle to the subgroup 
$\mathrm{G}$. The Levi-Civita connection is a $1$-form $Z$ on 
$\mathcal{F}(M^{n})$ with values in the Lie algebra $\mathg{so}(n)$.
We restrict the Levi-Civita connection to $\mathcal{R}$ and decompose 
it with respect to the decomposition of the Lie algebra $\mathg{so}(n)$,
\bdm
Z\big|_{T(\mathcal{R})} \ := \ Z^* \, \oplus \ \Gamma \, .
\edm
Then, $Z^*$ is a connection in the principal $\mathrm{G}$-bundle $\mathcal{R}$
and $\Gamma$ is a $1$-form on $M^{n}$ 
with values in the associated bundle $\mathcal{R} \times_{\mathrm{G}} 
\mathg{m}$. 
Suppose that the group $\mathrm{G}$ and the $\mathrm{G}$-structure are defined
by some differential form $\T$. Examples are
almost hermitian structures or almost metric contact structures. Then
the Riemannian covariant derivative of $\T$ is given by the formula
\bdm
\nabla^{\mathrm{LC}} \T \ = \ \vrho_*(\Gamma)(\T) \, ,
\edm
where $\vrho_*(\Gamma)(\T)$ denotes the algebraic action of the $2$-form $\Gamma$ on $\T$.
Some authors call $\Gamma$ the {\it intrinsic torsion}
of the $\mathrm{G}$-structure. There is a second notion, namely
the {\it characteristic connection} and the
{\it characteristic torsion} of a $\mathrm{G}$-structure. It is a
$\mathrm{G}$-connection $\nabla^{\mathrm{c}}$ with totally skew symmetric 
torsion tensor. Not
any type of geometric $\mathrm{G}$-structures admits a characteristic
connection. In order to formulate the condition, we embed the space of 
all $3$-forms into $\R^n \otimes \m$ using
the morphism
\bdm
\Theta \, : \, \Lambda^3(\R^n) \longrightarrow \R^n \otimes \m \, , \quad
\Theta(\T) \ := \ \sum_{i=1}^n e_i \otimes \mathrm{pr}_{\m}(e_i \haken
\T) \, .
\edm
A geometric $\mathrm{G}$-structure admits a characteristic connection 
$\nabla^{\mathrm{c}}$ if
and only if the intrinsic torsion $\Gamma$ belongs to the image of the 
$\Theta$. In this case, the intrinsic torsion is given by the equation (see \cite{Fri3}, \cite{Fri4})
\bdm
2 \, \Gamma \ = \ - \,  \Theta(\T^{\mathrm{c}}) \, .
\edm
For several geometric structures the characteristic torsion form 
has been computed 
explicitly in terms of the underlying geometric data. Formulas of that type
are known for almost hermitian structures, almost
metric contact structures and $\mathrm{G}_2$-structures in dimension $7$
(see \cite{FriedrichIvanov}). For a Riemannian naturally reductive
space $M^n = \mathrm{G}_1/\mathrm{G}$, we obtain a $\mathrm{G}$-reduction
$\mathcal{R} := \mathrm{G}_1\subset \mathcal{F}(M^{n})$ of the frame bundle.
Then the characteristic connection of the $\mathrm{G}$-structure
coincides with the {\it canonical} connection of the reductive space.
In this sense, we can understand the characteristic connection of a Riemannian 
$\mathrm{G}$-structure as a generalization of the canonical connection
of a Riemannian naturally reductive space. The canonical 
connection of a naturally reductive space has parallel torsion form and
parallel curvature tensor,
\bdm
\nabla^{\mathrm{c}} \T^{\mathrm{c}} \ = \ 0 \, , \quad
\nabla^{\mathrm{c}} \mathrm{R}^{\mathrm{c}} \ = \ 0 \, .  
\edm 
For arbitrary $\mathrm{G}$-structures and their characteristic connections,
these properties do not hold anymore. Corresponding examples are discussed
in \cite{FriedrichIvanov}. However, the parallelism of the torsion
form is an important property.
The first reason is that $\nabla^{\mathrm{c}} \T^{\mathrm{c}} = 0$
implies the conservation law $\delta(\T^ {\mathrm{c}}) = 0$, one of 
 the conditions for the NS-3-form in type II string theory (for constant 
dilaton). Moreover, if the torsion is parallel, several formulas for 
differential operators acting
on spinors simplify (see \cite{AgFr2}) and it is possible to investigate --
via integral formulas -- 
the space of parallel or harmonic spinors. Sasakian structures
or nearly K\"ahler structures 
have a parallel characteristic torsion form, even if they are not 
reductive. This motivates the investigation of Riemannian 
$\mathrm{G}$-structures with a parallel characteristic torsion form
in general.
In the present paper, we study the problem for almost hermitian
manifolds in dimension six.\\

\noindent
First we revisit almost hermitian manifolds in real
dimension six. The Hodge operator acts as a complex structure 
on $\Lambda^{3}(\R^6)$.
This observation simplifies, in dimension six, the description of 
the algebraic decomposition of the space of all $3$-forms 
$\Lambda^{3}(\R^6)$  and of the space  
$\R^6 \otimes \m^6$ containing the intrinsic torsion. We develop the 
algebraic part needed for the classification of
almost hermitian structures and we compute the
corresponding differential equations characterizing the  
sixteen classes of almost hermitian
manifolds (see \cite{GrayHervella}, \cite{FaFaSalamon} and 
\cite{ChiossiSalamon}). 
It is a basic property of $6$-dimensional
nearly K\"ahler manifolds that their characteristic torsion $\T^{\mathrm{c}}$
is $\nabla^{\mathrm{c}}$-parallel. The necessary formulas
proving that fact have been derived by the Japanese school at the beginning of the 70-ties of the last century (see \cite{Matsumoto}, \cite{Takamatsu},
\cite{YanoKon}). Later Gray (see \cite{Gray70}, \cite{Gray76}) and 
Kirichenko (see \cite{Kirichenko}) used these curvature identities 
for the investigation of the geometry of nearly K\"ahler manifolds. However,
the $\nabla^{\mathrm{c}}$-parallelism of the characteristic 
torsion $\T^{\mathrm{c}}$ was explicitly formulated only recently
(see \cite{BelgunMoroianu}, \cite{FriedrichIvanov}, \cite{Kirichenko}). 
We outline a short proof here, and continue our investigation 
along this path. Any almost hermitian manifold of 
type $\mathrm{G}_1$ admits a unique
characteristic connection (see \cite{FriedrichIvanov}).  
We study almost hermitian $\mathrm{G}_1$-manifolds with a parallel 
characteristic torsion. The $\U(3)$-orbit 
type of the characteristic torsion is constant. There are two possibilities.
If the vector part of the intrinsic torsion is non trivial, we obtain
two commuting Killing vector fields of constant length, and the
manifold is a torus fibration over some special $4$-manifold.
If the vector part vanishes, we list the relevant $\U(3)$-orbit types
of the torsion $3$-forms. It turns out that there exist only two
orbits with a non abelian isotropy (holonomy) group in dimension six.
These two types can be realized and the corresponding 
hermitian manifolds are twistor spaces or
the invariant, non K\"ahlerian hermitian structure on the Lie group 
$\SL(2,\C)$. Finally we classify all naturally reductive hermitian
$\mathcal{W}_3$-manifolds with  small isotropy group of the
characteristic torsion.
%
%
\section{Almost complex structures in real dimension four}\noindent
%
\subsection{The subgroup $\U(k)$ in $\SO(2k)$}\noindent
\vspace{3mm}

\noindent
We start with some notations that will be used throughout
this paper. $\R^n$  denotes the $n$-dimensional euclidian space. Using its 
scalar product $\langle \, , \, \rangle$, we
identify  euclidian space with its dual space, $\R^n \, = \, (\R^n)^*$.
$e_1 \, \ldots \, , e_n$ is an orthonormal basis in $\R^n$.
$\Lambda^l(\R^n)$ is the space of $l$-forms in $\R^n$.
$e_{i_1 \ldots i_l}$ means the exterior product $e_{i_1} \wedge \ldots \wedge
e_{i_l}$ 
of the corresponding $1$-forms.
We decompose a $2$-form $\omega$ or a $3$-form $\T$  
into their components,
\bdm
\omega \ = \ \sum_{1 \leq i < j \leq n} w_{ij} \cdot e_{ij} \, , \quad
\T \ = \ \sum_{1 \leq i < j < k \leq n} \T_{ijk} \cdot e_{ijk}  .
\edm 
The special orthogonal group $\SO(n)$ acts on $\Lambda^l(\R^n)$ and the
differential $\vrho_* : \mathg{so}(n) \rightarrow 
\mathrm{End}(\Lambda^l(\R^n))$ of this representation is given by
\bdm
\vrho_* (\omega) (\T) \ = \ \sum_{i=1}^n (e_i \haken \omega) \wedge
(e_i \haken \T) \ .
\edm
The space of $2$-forms
$\Lambda^2(\R^n) = \so(n)$ coincides with the Lie algebra of the special
orthogonal group, $\vrho$ is the 
adjoint representation and its differential $\vrho_*$ coincides with
the commutator action.\\
\vspace{2mm}

\noindent
We consider the complex structure $\J : \R^{2k} \rightarrow
\R^{2k}$ of the even-dimensional euclidian space. With respect to
the standard orthonormal basis it is given by
\bdm
\J e_{2i-1} \ = \ e_{2i} \, , \quad \J e_{2i} \ = \ - \, e_{2i-1} \, , \quad
i \ = \ 1 , \, 2 , \ldots  , \, k \ .
\edm
The subgroup $\U(k) \subset \SO(2k)$ consists of all orthogonal
transformations commuting with the complex structure
\bdm
\U(k) \ := \ \big\{A \in \SO(2k) \, : A \circ \J \, = \, \J \circ A \big\} \, .
\edm
The Lie algebra $\so(2k)$ splits into the subalgebra $\un(k)$ and
its orthogonal complement $\m$,
\bdm
\so(2k) \ = \ \Lambda^2(\R^{2k}) \ = \ \un(k) \, \oplus \m \, .
\edm
The complex structure $\J$ acts on $\Lambda^2(\R^{2k})$ 
as an involution. Using this involution, we can describe the spaces of the decomposition,
\bdm
\un(k) \ = \ \big\{ \omega \in \Lambda^2(\R^{2k}) \, : \, \J(\omega) \, = \, 
\omega \big\} \, , \quad 
\m \ = \  \big\{ \omega \in \Lambda^2(\R^{2k}) \, : \, \J(\omega) \, = \, 
- \, \omega \big\} \ . 
\edm
The center of the Lie algebra $\un(k)$ is generated by the $2$-form 
$\Omega(X,Y) := g(\J(X) \, , \, Y)$ and the Lie algebra splits into
\bdm
\un(k) \ = \ \su(k) \, \oplus \, \R^1 \cdot \Omega \, .
\edm
The Lie algebra $\un(k)$ is the space of all $2$-forms
defined by the equations,
\bdm
w_{2 i -1 , 2 j - 1} \, - \, w_{2 i , 2 j} \ = \ 0 \, , \quad - \, w_{2 i - 1
, 2 j} \, + \, w_{2 i, 2 j - 1} \ = \ 0 \, , \quad 1 \ \leq \ i \ < \ j \ \leq \
k \ .  
\edm 
The additional 
equation singling out the Lie algebra $\su(k)$ inside $\un(k)$ is
\bdm
w_{12} \, + \,  w_{34} \, + \, \ldots \,  + \, w_{2k-1,2k} \ = \ 0 \ .
\edm
%
%

\subsection{The decomposition of $\R^4 \otimes \m^2$}\noindent
\vspace{3mm}

\noindent
In dimension four, the Hodge operator as well as the complex structure
act on $2$-forms as involutions,
\bdm
\J^2 \ = \ \mathrm{Id} \ = \ *^2 \, , 
\quad \J \circ * \ = \ * \circ \J \, .
\edm
In contrast
to the higher-dimensional case, in real dimension four 
there are only two types. They are determined by the Nijenhuis tensor and
the differential of the K\"ahler form. 
In order to understand the geometric types of $\U(2)$-structures
on $4$-dimensional Riemannian manifolds, we need the decomposition of 
the representation $\R^4 \otimes \m^2$. Denote by
\bdm
\Phi \, : \, \R^4 \otimes \m^2 \longrightarrow \Lambda^3(\R^4) \, , \quad
\Phi(X \otimes \omega^2) \ := \ X \wedge \omega^2  
\edm
the total anti-symmetrization of a tensor in $\R^4 \otimes \m^2$. On the
other side, we embed the space of all $3$-forms into $\R^4 \otimes \m^2$ using
the morphism $\Theta  : \, \Lambda^3(\R^4) \rightarrow \R^4 \otimes \m^2$
defined in the introduction. A direct algebraic computation proves the 
following Lemma.
\begin{lem}
The morphism $\Phi : \R^4 \otimes \m^2 \rightarrow \Lambda^3(\R^4)$ 
is surjective and  $\Phi \circ \Theta$ acts on the space of all $3$-forms by
\bdm
\Phi \circ \Theta \ = \ \mathrm{Id} \ . 
\edm
\end{lem}
\noindent
Let us introduce two $\U(2)$-invariant subspaces of $\R^4 \otimes \m^2$,
\bdm
\mathcal{W}_2 \ := \ \mathrm{Ker}(\Phi) \, , \ 
\mathcal{W}_4 \ := \ \Theta\big(\Lambda^3(\R^4)\big) \ .
\edm
Obviously, $\R^4 \otimes \m^2$ splits under the action of the group
$\U(2)$ into these subspaces.
\begin{prop} 
$\mathcal{W}_2$ and $\mathcal{W}_4$ are real, irreducible $\U(2)$-representations.
\end{prop}
\begin{proof} We restrict the representation $\R^4 \otimes \m^2$ to the
subgroup $\SU(2)$. Then $\m^2$ is trivial and $\R^4 \otimes \m^2 = 
\R^4  \oplus  \R^4$ splits into two irreducible components under
the action of $\SU(2)$.
\end{proof}
%
%
\subsection{Geometric types of almost hermitian $4$-manifolds}\noindent
\vspace{3mm}

\noindent
Consider an almost hermitian manifold $(M^4, g, \J)$ and denote 
its Riemannian frame bundle by 
$\mathcal{F}(M^4)$. The almost hermitian structure is a reduction 
$\mathcal{R} \subset \mathcal{F}(M^{4})$ of the frame bundle to the subgroup 
$\U(2)$. The different non integrable types of hermitian structures 
are the irreducible 
components of the representation $\R^{4} \otimes \m^2$. We split the intrinsic
torsion $\Gamma$,
\bdm
\Gamma \ = \ \Gamma_4 \oplus \Gamma_4^* .
\edm
Note that, via the identification $\Theta$, $\Gamma_4$ is 
an ordinary $3$-form on the hermitian manifold. 
Moreover, in real dimension four, the differential and the co-differential
of the K\"ahler form coincide,
\bdm
\delta \, \Omega \ = \ - \, *  d \, *  \Omega \ = 
\ - \, *  d \, \Omega \, .
\edm
The co-differential of the K\"ahler form is given by the formula
\bdm
- \, \delta \, \Omega \ = \ \sum_{i=1}^4 e_i \haken \nabla^{\mathrm{LC}}_{e_i} 
\Omega \ = \ \sum_{i=1}^4 \Big\{ \Gamma(e_i)\big( e_i \haken
\Omega \, , \, - \big) \, - \, \Omega \big( e_i \haken \Gamma(e_i) \, , \, 
- \big) \Big\} \ =: \ \Pi(\Gamma) \, .
\edm
The map $\Pi : \R^4 \otimes \m^2 \rightarrow \R^4$ is obviously
$\U(2)$-equivariant. Consequently, the co-differential $\delta \, \Omega$
depends only on the $\Gamma_4$-part of the intrinsic torsion.
\begin{prop} 
Let $\Gamma_4 = \Theta(\T_4)$ be given by the $3$-form
$\T_4 \in \Lambda^3(\R^4)$. Then $\Pi(\Gamma_4) =  - \, 2 \,
\J(* \T_4)$ holds. In particular, the co-differential of the K\"ahler form 
of any almost hermitian $4$-manifold is given by the formula
\bdm
- \, *  d \, \Omega \ = \ \delta \, \Omega \ = \ 2 \, \J(* \T_4).
\edm
\end{prop}
\noindent
The Nijenhuis tensor in real dimension four has 
four components. Setting $\Ni_1 :=  \Ni(e_1 ,e_3)$, it may be written
in the form
\bdm
\Ni \ = \ (e_{13} \, - \, e_{24}) \otimes \Ni_1 \, - \, (e_{23} \, + \, 
e_{14}) \otimes \J(\Ni_1) \, .
\edm 
The anti-symmetrization map $\Phi : \R^4 \otimes
\m^2 \rightarrow \Lambda^3(\R^4)$ vanishes on the Nijenhuis tensor, i.e., 
$\Ni$ is an element of the subspace $\mathcal{W}_2$.
Consequently, there are two basic geometric types of almost hermitian $4$-manifolds.
They correspond to the Nijenhuis tensor (the $\Gamma_4^*$-part) and to the
differential $d \, \Omega$ of the K\"ahler form (the $\Gamma_4$-part). 
An almost hermitian
$4$-manifold admits a characteristic connection if and only if
its Nijenhuis tensor vanishes (hermitian manifold). In this case, the
characteristic torsion is given by the formula
$\T^{\mathrm{c}} =  - \, \J(d \, \Omega)$. It is 
$\nabla^{\mathrm{c}}$-parallel if
and only if the Lee form $\delta \Omega \circ \J$ is parallel with respect to
the
Levi-Civita connection. Hermitian manifolds of that type are called
generalized Hopf manifolds (see \cite{Vaisman}). The compact $4$-dimensional
generalized Hopf manifolds have been described by Belgun (see \cite{Belgun}).
%
\section{Almost complex structures in real dimension six}\noindent
In real dimension six, the Hodge operator as well as the complex structure
act on $3$-forms as complex structures. Moreover, the central element
$\Omega \in \un(3)$ acts on $3$-forms, too. These three operators split
the spaces $\Lambda^3(\R^6)$ and $\R^6 \otimes \m^6$ into $\U(3)$-irreducible
components. There are four basic types of almost hermitian $6$-manifolds.
They are characterized by the components of the derivative $d \, \Omega$ and
the Nijenhuis tensor $\Ni$.
%
\subsection{The decomposition of $\Lambda^3(\R^6)$}\noindent
%
\vspace{3mm}

\noindent
Two operators act on the space $\Lambda^3(\R^6)$, namely $\J$ and the Hodge 
operator $*$. The complex structure acts on a $3$-form $\T$ by
\bdm
(\J \,\T)(X \, ,\, Y \,,\, Z)\ :=\ \T(\J X \, ,\,\J Y\, ,\,\J Z)\, .
\edm 
We obtain a $(\Z_4 \oplus \Z_4)$-action on the space of all
$3$-forms,
\bdm
\J^2 \ = \ - \, \mathrm{Id} \, , \quad *^2 \ = \ - \, \mathrm{Id} \, , 
\quad \J \circ * \ = \ * \circ \J \, .
\edm
Let us decompose $\Lambda^3(\R^6)$ into two $\U(3)$-invariant subspaces,
\bdm
\Lambda^3(\R^6) \ := \ \Lambda^3_+(\R^6) \, \oplus \,\Lambda^3_-(\R^6),
\quad \text{where} \quad  \Lambda^3_{\pm} \ := \ \big\{ \T \in 
\Lambda^3(\R^6) \, : \, \J(\T) \, = \, \pm \, * \T \big\} \, .
\edm
We embed the standard representation $\R^6$ into the $3$-forms,
\bdm
\Lambda^3_6(\R^6) \ := \ \big\{ X \wedge \Omega \, : \, X \in 
\R^6 \big\} \, .
\edm
\begin{lem} 
The spaces $\Lambda^3_-(\R^6) = \Lambda^3_6(\R^6)$ coincide. 
The space $\Lambda^3_+(\R^6)$ is its orthogonal complement,
\bdm
\Lambda^3_+(\R^6) \ = \ \big\{ \T \, : \, \Omega \wedge \T = 0 \big\} \ = \ 
\big\{ \T \, : \, \Omega \wedge * \T = 0 \big\} \ = \ \big\{ \T \, : \, 
\J(\T) =  * \T \big\} \ .
\edm
\end{lem}
\noindent
The $\U(3)$-representation $\Lambda^3_+(\R^6)$ is not irreducible. In order
to decompose it, we consider the action of the central element 
$\Omega \in \un(3)$ on the space $\Lambda^3(\R^6)$. We will denote 
by $\tau : \Lambda^3(\R^6) \rightarrow \Lambda^3(\R^6)$ this special
antisymmetric operator acting on $3$-forms, 
\bdm
\tau(\T) \ := \ \vrho_* (\Omega) (\T) \ = \ \sum_{i=1}^6 (e_i \haken \Omega) 
\wedge (e_i \haken \T) \, .
\edm
\begin{lem} 
The symmetric endomorphism $\tau^2$ has two eigenvalues and splits
the space $\Lambda^3_+(\R^3)$ into a $2$-dimensional and a $12$-dimensional 
space
\bdm
\Lambda^3_+(\R^6) \ := \ \Lambda^3_2(\R^6) \oplus \Lambda^3_{12}(\R^6),
\edm
where
\begin{eqnarray*} 
\Lambda^3_2(\R^6) &:=& \big\{\T \in \Lambda^3(\R^6) :  
\tau^2(\T) \, = \, - \, 9 \, \T \, , \ \J(\T) \, = \, * \T \big\} \\
\Lambda^3_{12}(\R^6) &:=&  \big\{\T \in \Lambda^3(\R^6) :  
\tau^2(\T) \, = \, - \, \T \, , \ \J(\T) \, = \, * \T \big\} \ .
\end{eqnarray*}
The antisymmetric endomorphism $\tau$ preserves any of these spaces
and acts as
\bdm
\tau(\T) \, = \, 3 \, * \T \ \, \mathrm{on} \ \, \Lambda^3_2(\R^6) \, , \quad
\tau(\T) \, = \, * \T \ \,  \mathrm{on} \ \, \Lambda^3_6(\R^6) \, , \quad
\tau(\T) \, = \, - \, * \T \ \, \mathrm{on} \ \, \Lambda^3_{12}(\R^6) \, .
\edm
\end{lem}
\noindent
It is useful to have at hand an explicit basis in any of these spaces: 
\bdm
\hspace{0.3cm} \mathrm{in} \  \Lambda^3_2(\R^6) : \quad
- \,e_{246} \, + \, e_{136} \, + \, e_{145} \, + \, e_{235} \, , \quad
- \,e_{135} \, + \, e_{245} \, + \, e_{236} \, + \, e_{146} \ . 
\edm
\begin{eqnarray*}
\hspace{-3cm}
\mathrm{in} \ \Lambda^3_6(\R^6) : \quad
&& e_{134} \, + \, e_{156} \, , \quad  e_{234} \, + \, e_{256} \, , 
\quad e_{123} \, + \, e_{356} \, , \\  
&& e_{124} \, + \, e_{456} \, ,
\quad e_{125} \, + \, e_{345} \, , \quad  e_{126} \, + \, e_{346} \, . 
\end{eqnarray*}
\vspace{0.5mm}
  
\begin{eqnarray*}
\mathrm{in} \ \Lambda^3_{12}(\R^6) : \quad
&& e_{123} \, - \, e_{356} \, , \quad  e_{124} \, - \, e_{456} \, , 
\quad e_{125} \, - \, e_{345} \, , \quad  e_{126} \, - \, e_{346} \, , \\
&& e_{134} \, - \, e_{156} \, , \quad 
e_{234} \, - \, e_{256} \, , \quad
e_{135} \, + \, e_{245} \, , \quad  e_{246} \, + \, e_{136} \, , \\ 
&& e_{135} \, + \, e_{236} \, , \quad  e_{246} \, + \, e_{145} \, , 
\quad e_{135} \, + \, e_{146} \, , 
\quad e_{246} \, + \, e_{235} \ . 
\end{eqnarray*}
\noindent
Any $2$-dimensional real representation of the simply connected 
compact Lie group $\SU(3)$ is trivial. Therefore, the subgroup $\SU(3)$
preserves any $3$-form in $\Lambda^3_2(\R^6)$ and we
can understand the subgroup $\SU(3) \subset \U(3)$ as the isotropy
group of a $3$-form of that type. Moreover, we obtain
a $\SU(3)$-isomorphism between $\R^6$ and $\m^6$.
\begin{cor} 
For any non trivial $3$-form $\T \in \Lambda^3_2(\R^6)$,
the map $X \rightarrow X \haken \T$ defines an $\SU(3)$-isomorphism
between $\R^6$ and $\m^6$. 
\end{cor} 
\begin{NB}
The irreducible $6$-dimensional $\U(3)$-representation
 $\m^6$ is not equivalent to the standard representation in $\R^6$. Indeed,
$\J$ is an element of the group $\U(3)$ and we compute its trace 
in $\R^6$ and in $\m^6$, $\mathrm{Tr}_{\R^6}(\J) = 0 \, , \
\mathrm{Tr}_{\m^6}(\J) =  - \, 6 $.
\end{NB}
\noindent
Now we prove that $\Lambda^3_{12}(\R^6)$ is 
irreducible. We use the fact that there exists only one non trivial 
$6$-dimensional real representation of the group $\SU(3)$. For completeness, 
we sketch the proof, too. 
\begin{prop} \label{6-dimensionaleDarstellung} 
Any $6$-dimensional real representation of the group
$\SU(3)$ is either trivial or isomorphic to the standard representation
in $\R^6 = \C^3$.
\end{prop}
\begin{proof} 
The euclidian group $\SO(5)$ does not contain a subgroup of
dimension eight. Consequently, any $6$-dimensional
real representation $V^6$ of $\SU(3)$ is either trivial or irreducible.
Suppose that $V^6 \neq \R^6$ is irreducible. Since $\C^3$ and its 
complex conjugation are the only complex irreducible
representation of the group $\SU(3)$ in dimension three, the complexification
$(V^6)^{\C}$ must be irreducible (see \cite[Lemma 3.58]{Adams}). Again, there
are only two $6$-dimensional irreducible $\SU(3)$-representations, namely
$\mathrm{Sym}^2(\C^3)$ and its conjugation. We compute
the character of the element $g = \mathrm{diag}(z,z,z^{-2}) \in \SU(3)$,
\bdm
\chi_{\mathrm{Sym}^2(\C^3)}(g) \ = \ 3 \, z^2 \, + \, z^{-4} \, + \, 2 \, 
z^{-1} \ ,
\edm
and conclude that $\mathrm{Sym}^2(\C^3)$ is not a real representation.  
\end{proof}
\begin{thm} The decomposition
\bdm
\Lambda^3(\R^6) \ = \ \Lambda^3_2(\R^6) \oplus \Lambda^3_6(\R^6) \oplus
\Lambda^3_{12}(\R^6)
\edm
splits the space of all $3$-forms into 
irreducible, real $\U(3)$-representations. Moreover, $\Lambda^3_6(\R^6)$ and 
$\Lambda^3_{12}(\R^6)$ are irreducible $\SU(3)$-representations. 
$\Lambda^3_{2}(\R^6)$ is the trivial $2$-dimensional real 
$\SU(3)$-representation. $\Lambda^3_2(\R^6)$, $\Lambda^3_6(\R^6)$  
and  $\Lambda^3_{12}(\R^6)$ are irreducible, complex representations 
of dimensions $1$, $3$ and $6$, respectively.  
\end{thm}
\begin{proof} 
The $\SU(3)$-representation $\Lambda^3_{12}(\R^6)$ can split 
only into $\R^6 \, \oplus \, \R^6$ or $\R^6 \, \oplus \, 6 \, \R^1$ 
(see Proposition \ref{6-dimensionaleDarstellung}). Consider the following 
elements of the group $\SU(3)$,
\bdm
g_1 \ := \
\left[\ba{cccccc} -1 & 0 & 0 & 0 & 0 & 0 \\ 
0 & -1 & 0 & 0 & 0 & 0 \\  
0 & 0 & -1 & 0 & 0 & 0 \\
0 & 0 & 0 & -1  & 0 & 0 \\
0 & 0 & 0 & 0 & 1 & 0 \\
0 & 0 & 0 & 0 & 0 & 1 \ea\right] \ , \quad
g_2 \ := \
\left[\ba{cccccc} 0 & 0 & 0 & 1 & 0 & 0 \\ 
0 & 0 & -1 & 0 & 0 & 0 \\  
0 & 1 & 0 & 0 & 0 & 0 \\
-1 & 0 & 0 & 0  & 0 & 0 \\
0 & 0 & 0 & 0 & 1 & 0 \\
0 & 0 & 0 & 0 & 0 & 1 \ea\right]  .
\edm
We compute the values of the characters,
\begin{eqnarray*} 
\chi_{\Lambda^3_{12}(\R^6)}(g_1) &=& 4 \, , \quad
\chi_{\R^6 \oplus \R^6}(g_1) \ = \ - \, 4 \, , \\
\chi_{\Lambda^3_{12}(\R^6)}(g_2) &=& 0 \, , \quad
\chi_{\R^6 \oplus 6 \, \R^1}(g_2) \ = \ 8 \ ,
\end{eqnarray*}
i.e., both cases are impossible.
\end{proof} 
\subsection{The decomposition of $\R^6 \otimes \m^6$}\noindent
\vspace{3mm}

\noindent
In order to understand the geometric types of $\U(3)$-structures
on $6$-dimensional Riemannian manifolds, we need the decomposition of 
the representation $\R^6 \otimes \m^6$. Denote by
\bdm
\Phi \, : \, \R^6 \otimes \m^6 \longrightarrow \Lambda^3(\R^6) \, , \quad
\Phi(X \otimes \omega^2) \ := \ X \wedge \omega^2  
\edm
the total anti-symmetrization of a tensor in $\R^6 \otimes \m^6$. On the
other side, we embed the space of all $3$-forms into $\R^6 \otimes \m^6$ using
the morphism
\bdm
\Theta \, : \, \Lambda^3(\R^6) \longrightarrow \R^6 \otimes \m^6 \, , \quad
\Theta(\T) \ := \ \sum_{i=1}^6 e_i \otimes \mathrm{pr}_{\m^6}(e_i \haken
\T) \, .
\edm
A direct algebraic computation yields the following Lemma.
\begin{lem}
The morphism $\Phi : \R^6 \otimes \m^6 \rightarrow \Lambda^3(\R^6)$ 
is surjective and  $\Phi \circ \Theta$ acts on the space of all $3$-forms by
\bdm
\Phi \circ \Theta \ = \ 3 \, \mathrm{Id} \, \ \mathrm{on} \, \
\Lambda^3_2(\R^6) \, , \quad 
\Phi \circ \Theta \ = \ \mathrm{Id} \, \ \mathrm{on} \, \
\Lambda^3_6(\R^6) \oplus \Lambda^3_{12}(\R^6) \ . 
\edm
\end{lem}
\noindent
Let us introduce four $\U(3)$-invariant subspaces of $\R^6 \otimes \m^6$,
\bdm
\mathcal{W}_1 \ := \ \Theta\big(\Lambda^3_2(\R^6)\big) \, , \
\mathcal{W}_2 \ := \ \mathrm{Ker}(\Phi) \, , \
\mathcal{W}_3 \ := \ \Theta\big(\Lambda^3_{12}(\R^6)\big) \, , \
\mathcal{W}_4 \ := \ \Theta\big(\Lambda^3_6(\R^6)\big) \ .
\edm
$\R^6 \otimes \m^6$ splits under the action of the group
$\U(3)$ into these subspaces. We investigate the representation
$\mathcal{W}_2$. It splits as a $\SU(3)$-representation. Fix a
$3$-form $\T \in \Lambda^3_2(\R^6)$. The group $\SU(3)$ stabilizes
$\T$ and the morphism
\bdm
\Psi_{\T} \, : \, \R^6 \otimes \m^6 \longrightarrow \Lambda^2(\R^6) \, ,
\quad \Psi_{\T}(X \otimes \omega^2) \ := \ *\big( (X \haken \T) 
\wedge \omega^2) \big) 
\edm 
is $\SU(3)$-equivariant. We can control the image of $\Psi_{\T}$.
\begin{lem} 
For any non trivial form $\T \in \Lambda^3_2(\R^6)$, 
the image of $\Psi_{\T}$ is contained in the Lie algebra 
$\un(3)$. Moreover, $\Psi_{\T}$ maps $\mathcal{W}_2$ surjectively
onto the Lie algebra $\su(3)$.
\end{lem}
\noindent
Now we decompose the representation $\mathcal{W}_2$ under
the action of the group $\SU(3)$.
\begin{thm} 
Fix two linearly independent $3$-forms $\T_1 \, , \, 
\T_2$ in $\Lambda^3_2(\R^6)$. Then the map
\bdm
\Psi_{\T_1} \oplus \Psi_{\T_2} \ : \ \mathcal{W}_2 \longrightarrow 
\su(3) \oplus \su(3)
\edm
is an isomorphism of  $\SU(3)$-representations.
\end{thm}
\noindent
Finally we prove that $\mathcal{W}_2$ is $\U(3)$-irreducible.
\begin{thm} 
$\mathcal{W}_2$ is a real, irreducible $\U(3)$-representation
of dimension $16$.
\end{thm}
\begin{proof} 
The element $e_1  \otimes  (e_{14}  +  e_{23})  +  
e_2  \otimes  (e_{13}  -  e_{24})  \in  \mathcal{W}_2$ is not
invariant under the action of the $1$-parameter group generated by
the central element $\Omega \in \un(3)$. Suppose that  $\mathcal{W}_2$ 
is $\U(3)$-reducible. Then the adjoint representation of $\SU(3)$
extends to a representation $\kappa$ of the group $\U(3)$. In particular, 
$\Omega \in \un(3)$ defines a non trivial, skew symmetric $\SU(3)$-invariant
operator $\kappa_*(\Omega) : \su(3) \rightarrow \su(3)$. Since for
any simple Lie group $\mathrm{G}$ we have
\bdm
\Lambda^2(\g)^{\mathrm{G}} \ = \ 0 \, ,
\edm  
this is a contradiction.
\end{proof}
\begin{cor} 
The $\U(3)$-representation $\R^6 \otimes \m^6$ splits
into four irreducible representations,
\bdm
\R^6 \otimes \m^6 \ = \ \mathcal{W}_1 \oplus \mathcal{W}_2 \oplus 
\mathcal{W}_3 \oplus \mathcal{W}_4 .
\edm
The  $\SU(3)$-representation $\R^6 \otimes \m^6$ splits into
\bdm
\R^6 \otimes \m^6 \ = \ \R^2 \oplus \big(\su(3) \oplus \su(3)\big)
\oplus \mathcal{W}_3 \oplus \mathcal{W}_4.
\edm
\end{cor}
%
\subsection{The sixteen classes of almost hermitian structures}\noindent
\vspace{3mm}

\noindent
Consider an almost hermitian manifold $(M^6, g, \J)$ and denote 
its Riemannian frame bundle by 
$\mathcal{F}(M^6)$. The almost hermitian structure is a reduction 
$\mathcal{R} \subset \mathcal{F}(M^{6})$ of the frame bundle to the subgroup 
$\U(3)$. We restrict the Levi-Civita connection to $\mathcal{R}$ and decompose 
it with respect to the decomposition of the Lie algebra $\mathg{so}(6)$:
\bdm
Z\big|_{T(\mathcal{R})} \ := \ Z^* \, \oplus \, \Gamma \, .
\edm
The
Riemannian covariant derivative of the K\"ahler form is given by the formula
\bdm
\big(\nabla^{\mathrm{LC}}_X \Omega\big)\big(Y\, , \, Z \big) \ = \ 
\Gamma(X)\big( Y \haken \Omega \, , \, Z \big) \, - 
\, \Omega \big( Y \haken \Gamma(X) \, , \, Z \big) \, .
\edm
The basic types of hermitian structures are the irreducible 
components of the representation $\R^{6} \otimes \m^6$. We split the
intrinsic torsion $\Gamma$,
\bdm
\Gamma \ = \ \Gamma_2 \oplus \Gamma_6 \oplus \Gamma_{12} \oplus \Gamma_{16} .
\edm
Note that, via the identification $\Theta$, $\Gamma_2$ and $\Gamma_{12}$ 
are $3$-forms on 
the hermitian manifold and $\Gamma_6 = \Theta(X \wedge \Omega)$
is a vector field.
%
\subsection{The co-differential $\delta \, \Omega$}\noindent
\vspace{3mm}

\noindent
The co-differential of any exterior form $\alpha$ on a Riemannian manifold
is given by the formula
\bdm
\delta\, \alpha \ = \ - \, \sum_{i=1}^n e_i \haken \nabla^{\mathrm{LC}}_{e_i}
\alpha \, .
\edm
Inserting the formula for the covariant derivative of the K\"ahler form, we
obtain
\bdm
- \, \delta \, \Omega \ = \ \sum_{i=1}^6 \Big\{ \Gamma(e_i)\big( e_i \haken
\Omega \, , \, - \big) \, - \, \Omega \big( e_i \haken \Gamma(e_i) \, , \, 
- \big) \Big\} \ =: \ \Pi(\Gamma) \, .
\edm
The map $\Pi : \R^6 \otimes \m^6 \rightarrow \R^6$ is obviously
$\U(3)$-equivariant. Consequently, the co-differential $\delta \, \Omega$
depends only on the $\mathcal{W}_4$-part of the intrinsic torsion. We 
compute the relation explicitly.
\begin{prop} 
Let $\Gamma_6 = \Theta(X \wedge \Omega)$ be given by the vector
$X \in \R^6$. Then $\Pi(\Gamma_6) =  4 \, X$ holds. In particular, 
the co-differential of the K\"ahler form of any almost hermitian manifold
is given by the formula
\bdm
\delta \, \Omega \ = \ - \, 4 \, X \, \quad \mathrm{where} \quad
\Gamma_6 \ = \ \Theta(X \wedge \Omega) \, .
\edm
\end{prop}
\subsection{The differential $d \, \Omega$}\noindent
\vspace{3mm}

\noindent
We handle the differential of the K\"ahler form in a similar way. Indeed,
the differential of an arbitrary exterior form on a Riemannian manifold
can be computed by the formula
\bdm
d \, \alpha \ = \ \sum_{i=1}^n e_i \wedge \nabla^{\mathrm{LC}}_{e_i}\alpha\, .
\edm
Inserting again the formula for the covariant derivative of the K\"ahler form, 
we obtain
\bdm
d \, \Omega \ = \ \frac{1}{2}\sum_{i,j=1}^6 e_i \wedge e_j \wedge 
\Big\{ \Gamma(e_i)\big( e_j \haken \Omega \, , \, - \big) \, - \, \Omega 
\big( e_j \haken \Gamma(e_i) \, , \, 
- \big) \Big\} \ =: \ \Pi_1(\Gamma) \, .
\edm
The map $\Pi_1 : \R^6 \otimes \m^6 \rightarrow \Lambda^3(\R^6)$ is obviously
$\U(3)$-equivariant. Consequently, the differential $\d \, \Omega$
depends only on the $\Theta(\Lambda^3(\R^6))$-part of the intrinsic torsion.
Moreover, we need a formula for the endomorphism $\Pi_1 \circ \Theta : 
\Lambda^3(\R^6) \rightarrow \Lambda^3(\R^6)$.
\begin{prop} 
The endomorphism $\Pi_1 \circ \Theta$ is given on the 
irreducible components by the formulas
\begin{enumerate}
\item $\Pi_1 \circ \Theta (\T) = - \, 6 \, * \T \ $ for $\T \in 
\Lambda^3_2(\R^6)$.
\item $\Pi_1 \circ \Theta (\T) = - \, 2 \, * \T \ $ for $\T \in 
\Lambda^3_6(\R^6)$.
\item $\Pi_1 \circ \Theta (\T) = \ \ \, 2 \, * \T \ $ for $\T \in 
\Lambda^3_{12}(\R^6)$.
\end{enumerate}
\end{prop} 
\vspace{3mm}

\noindent
Let us summarize the result of these algebraic computations.
\begin{thm}\label{AbleitungOmega} 
Let $(M^6, g, \J)$ be an almost hermitian manifold of type
\bdm
\Gamma \ = \ \Gamma_2 \oplus \Gamma_6 \oplus \Gamma_{12} \oplus \Gamma_{16} .
\edm
Suppose that the first three parts of the intrinsic torsion are given
$3$-forms in the corresponding component of $\Lambda^3(\R^6)$,
\bdm
\Gamma_2 \ = \ \Theta(\T_2) \, , \quad \Gamma_6 \ = \ \Theta(X \wedge \Omega)
\, , \quad \Gamma_{12} \ = \ \Theta(\T_{12}).
\edm
The differential and the co-differential of the K\"ahler form do not
depend on the $\mathcal{W}_2$-component of the intrinsic torsion.
Moreover, we have
\bdm
\delta \, \Omega \ = \ - \, 4 \, X \, , \quad
d \, \Omega \ = \ - \, 6 \, * \T_2 \, - \, 2 \, * (X \wedge \Omega) \, 
+ \, 2 \, * \T_{12}.
\edm
\end{thm}

\subsection{The Nijenhuis tensor}\noindent
\vspace{3mm}

\noindent
The Nijenhuis tensor 
\bdm
\Ni(X \, , \, Y) \ := \ \big[ \J(X) \, , \, \J(Y) \big] \, - \, 
\J \big[ X \, , \, \J(Y) \big] \, - \, \J \big[ \J(X) \, , \, Y \big]
\, - \, \big[ X \, , \, Y \big]
\edm
in real dimension six has $18$ components, 
\bdm
\Ni_1 \ := \ \Ni(e_1 \, , \, e_3) \, , \quad
\Ni_2 \ := \ \Ni(e_1 \, , \, e_5) \, , \quad
\Ni_3 \ := \ \Ni(e_3 \, , \, e_5) \, , \quad
\edm
and is given by
\begin{eqnarray*}
\Ni &=&(e_{13} \, - \, e_{24}) \otimes \Ni_1 \, - \, (e_{23} \, + \, 
e_{14}) \otimes \J(\Ni_1) \, + \, 
(e_{15} \, - \, e_{26}) \otimes \Ni_2 \\
&-&(e_{25} \, + \, 
e_{16}) \otimes \J(\Ni_2) \, + \, 
(e_{35} \, - \, e_{46}) \otimes \Ni_3 \, - \, (e_{36} \, + \, 
e_{45}) \otimes \J(\Ni_3) . 
\end{eqnarray*} 
We apply the anti-symmetrization map $\Phi : \R^6 \otimes
\m^6 \rightarrow \Lambda^3(\R^6)$. Then $\Phi(\Ni)$ is contained in 
$\Lambda^3_2(\R^6)$. Consequently, the Nijenhuis tensor is an element
of the subspace $\mathcal{W}_1 \, \oplus \mathcal{W}_2 \ \subset \R^6 \otimes 
\m^6$ and coincides with the $(\mathcal{W}_1 \oplus 
\mathcal{W}_2)$-part of the intrinsic torsion. 
In particular, we obtain a characterization of hermitian manifolds.
\begin{thm} 
The almost complex structure $\J$ 
is integrable if and only if the $(\mathcal{W}_1 \oplus \mathcal{W}_2)$-part 
of its
intrinsic torsion vanishes. The Nijenhuis tensor is totally skew symmetric
if and only if the $\mathcal{W}_2$-part of the intrinsic torsion vanishes.
\end{thm}
%
\subsection{Differential equations characterizing the types}\noindent
\vspace{3mm}

\noindent
We identified the different parts of the intrinsic torsion with
the differential and the co-differential of the K\"ahler form as well as
with the Nijenhuis tensor. These formulas yield differential equations
characterizing any type of a non integrable hermitian geometry. Some
of these classes have special names. In general, we fix a $6$-dimensional
almost hermitian manifold $(M^6, g , \J)$.
\begin{cor} The following conditions are equivalent:
\begin{enumerate}
\item The structure is of type $\mathcal{W}_1 \oplus \mathcal{W}_2 
\oplus \mathcal{W}_3$.
\vspace{1mm}

\item $\delta \, \Omega = 0$.
\vspace{1mm}

\item $\Omega \wedge d \, \Omega = 0$.
\vspace{1mm}

\item $\J(d \, \Omega) =  * d\, \Omega$.
\end{enumerate}
\end{cor}
\noindent
Manifolds of that type are called {\it almost
semi-K\"ahler} or {\it co-symplectic}.
\begin{cor} 
The following conditions are equivalent:
\begin{enumerate}
\item The structure is of type $\mathcal{W}_1 \oplus \mathcal{W}_2 
\oplus \mathcal{W}_4$.
\vspace{1mm}

\item $\tau^2\big[d \, \Omega - \frac{1}{2} \, *(\delta \, \Omega
\wedge \Omega) \big] = - \, 9 \, \big[d \, \Omega - \frac{1}{2} \, 
*(\delta \, \Omega \wedge \Omega) \big] $.
\end{enumerate}
\end{cor}
\begin{cor}\label{W134} 
The following conditions are equivalent:
\begin{enumerate}
\item The structure is of type $\mathcal{W}_1 \oplus \mathcal{W}_3 
\oplus \mathcal{W}_4$.
\vspace{1mm}

\item The Nijenhuis tensor is totally skew symmetric.
\vspace{1mm}

\item There exists a linear connection $\nabla$ preserving the
almost hermitian structure and with totally skew symmetric torsion.
\end{enumerate}
\end{cor}
\noindent
The equivalence of the second and third condition
has been proved in \cite{FriedrichIvanov}, see also \cite{Fri3}.
Almost hermitian manifolds satisfying the latter condition 
are called {\it $\mathrm{G}_1$-manifolds}.
\begin{cor} 
The following conditions are equivalent:
\begin{enumerate}
\item The structure is of type $\mathcal{W}_2 \oplus \mathcal{W}_3 
\oplus \mathcal{W}_4$.
\vspace{1mm}

\item $\tau^2(d \, \Omega) = - \, d \, \Omega$.
\end{enumerate}
\end{cor}
\noindent
Almost hermitian manifolds of that type are called 
{\it $\mathrm{G}_2$-manifolds}.\\

\noindent
We investigate next the almost hermitian structures of pure type, where
only one component of the intrinsic torsion does not vanish.
\begin{cor} The following conditions are equivalent:
\begin{enumerate}
\item The structure is of type $\mathcal{W}_2$.
\vspace{1mm}

\item $d \, \Omega = 0$.
\end{enumerate}
\end{cor}
\noindent
Manifolds of that type are called {\it almost K\"ahler} or 
{\it symplectic}.
\begin{cor} The following conditions are equivalent:
\begin{enumerate}
\item The structure is of type $\mathcal{W}_3$.
\vspace{1mm}

\item $\J$ is integrable and $\delta \, \Omega = 0$.
\vspace{1mm}

\item $\Ni$ is totally skew symmetric and $\J(d \, \Omega) = \, * d\, \Omega \, , \,
\tau^2(d \, \Omega) = - \, d \, \Omega$.
\end{enumerate}
\end{cor}
\noindent
Almost hermitian manifolds of that type are called {\it semi-K\"ahler}.
\begin{cor} The following conditions are equivalent:
\begin{enumerate}
\item The structure is of type $\mathcal{W}_4$.
\vspace{1mm}

\item $\Ni$ is totally skew symmetric and $2 \, d \, \Omega = (\delta \, 
\Omega \circ \J) \wedge \Omega$.
\end{enumerate}
\end{cor}
\noindent
Manifolds of that type are called {\it locally conformal K\"ahler}.\\

\noindent
The most interesting and rigid class of almost hermitian manifolds in dimension
six is the class of so called {\it nearly K\"ahler manifolds}. In the sixties
and seventies of the last century, they have also been called
{\it Tachibana spaces} or {\it $\mathrm{K}$-spaces} (see \cite{Gray70}, 
\cite{Takamatsu}, \cite{Matsumoto}, \cite{YanoKon}). Nearly K\"ahler 
manifolds correspond to the pure type
$\mathcal{W}_1$ and we describe this class of almost
hermitian structures in the spirit of the previous corollaries. 
\begin{cor} 
The following conditions are equivalent:
\begin{enumerate}
\item The structure is of type $\mathcal{W}_1$.
\vspace{1mm}

\item $\Ni$ is totally skew symmetric and $\delta \, \Omega = 0 \, , \ 
\tau^2(d \, \Omega) = - \, 9 \, d \, \Omega$.
\end{enumerate}
\end{cor}
\noindent
Furthermore, the differential of the K\"ahler form satisfies the following
equations,
\bdm
\Omega \wedge d \, \Omega \ = \ 0 \, , \quad \J(d \, \Omega) \ = \  
* d \, \Omega  .
\edm
\noindent
There is an equivalent characterization of nearly K\"ahler manifolds.
\begin{thm}
An almost hermitian manifold is nearly K\"ahler if and only if, for
any vector $X$,
\bdm
\big(\nabla^{\mathrm{LC}}_X \Omega \big)\big(X \, , \, - \big) \ = \ 0 \, .
\edm
\end{thm}
\begin{proof} 
Consider the $\U(3)$-equivariant 
map $\ \R^6 \otimes \m^6 \rightarrow 
\mathrm{S}^2(\R^6) \otimes \R^6 \ $ defined by
\bdm
\Gamma \longrightarrow \hat{\Gamma}(X \, , \, Y) \ = \ 
\big(\nabla^{\mathrm{LC}}_X \Omega \big)\big(Y\big) \, + \, 
\big(\nabla^{\mathrm{LC}}_Y \Omega \big)\big(X\big) \, . 
\edm
It turns out that its kernel coincides with the subspace $\mathcal{W}_1
\subset \R^6 \otimes \m^6$.
\end{proof}
%
%
%
\section{The characteristic connection of a $\mathrm{G}_1$-manifold}\noindent
\subsection{The general formula for the characteristic connection}\noindent
\vspace{3mm}

\noindent
An almost hermitian manifold is of type $\mathrm{G}_1$ if and only if 
it admits a linear connection preserving the
structure  with skew symmetric torsion (see Corollary \ref{W134}), and,
in this case, the connection is unique. In this 
generality, this result has been proved in the paper \cite{FriedrichIvanov}. 
For special types of almost hermitian manifolds, the {\it characteristic 
connection} has been considered before. For nearly K\"ahler manifolds,
 A.~Gray used it in 1970 (see \cite[page 304]{Gray70}) in order to express
the Chern classes. In 1976 (see \cite[page 237]{Gray76}),
he proved that the first Chern class of a $6$-dimensional nearly
K\"ahler, non-K\"ahler manifold vanishes. On the other hand,
the characteristic connection of a hermitian manifold 
has been used by J.-M.~Bismut in 1989
in the proof of the local index theorem (see \cite{Bismut}). Let us
compute the formula for the torsion of the characteristic connection
of an almost hermitian manifold of type $\mathrm{G}_1$. Using the ansatz
\bdm
\Gamma \ = \ \Gamma_2  \oplus  \Gamma_6  \oplus  \Gamma_{12} \, , 
\quad \Gamma_2 \ = \ \Theta(\T_2) \, , \quad \Gamma_6 \ = \ \Theta(X \wedge
\Omega) \, , \quad \Gamma_{12} \ = \ \Theta(\T_{12}) 
\edm
as well as the formula $2 \, \Gamma =  - \, \Theta(\T^{\mathrm{c}})$ 
relating $\Gamma$ and the torsion form of the characteristic connection 
$\nabla^{\mathrm{c}}$
(see \cite{Fri3}, \cite{Fri4}), we obtain by Theorem \ref{AbleitungOmega} 
\bdm
\T^{\mathrm{c}} \ = \ - \, 2 \, \T_2 \, - \, 2 \, 
(X \wedge \Omega) \, - \, 2 \, \T_{12} 
 \ = \ - \, 8 \, \T_2 \, + \, \J(d \, \Omega)  .
\edm
The torsion form of the characteristic connection of a hermitian manifold 
($\Gamma_2 = \Gamma_{16} = 0)$ is the twisted
differential of the K\"ahler form (see \cite[Theorem $10.1$]{FriedrichIvanov}),
\bdm
\T^{\mathrm{c}} \ = \  \J(d \, \Omega) .
\edm
For nearly K\"ahler manifolds, we have
$d \, \Omega = - \, 6 \, * \T_2 = - \, 6 \, \J(\T_2)$. The torsion
of the characteristic connection is again proportional to the twisted
differential of the K\"ahler form,
\bdm
\T^{\mathrm{c}} \ = \ - \, \frac{1}{3} \, \J(d \, \Omega )  .
\edm
An easy computation yields an equivalent formula for the characteristic
connection and its torsion, namely
\bdm
\T^{\mathrm{c}}(X \, , \, Y) \ = \ - \, \J\big((\nabla^{\mathrm{LC}}_X \J)(Y)\big) \, , \quad
\nabla^{\mathrm{c}}_X Y \ = \ \frac{1}{2} \, \big( 
\nabla^{\mathrm{LC}}_X Y \, - 
\, \J(\nabla^{\mathrm{LC}}_X \J(Y)) \big) .
\edm
The latter formula is the original definition of the characteristic 
connection of a nearly K\"ahler manifold as it appears in the papers of 
A.~Gray (\cite{Gray70}, \cite{Gray76}).\\

\noindent
Combining the formula
\bdm
\T^{\mathrm{c}} \ = \ - \, 2 \, \T_2 \, - \, 2 \, 
(X \wedge \Omega) \, - \, 2 \, \T_{12} 
\edm
with the general formula of Theorem \ref{AbleitungOmega} 
\bdm
d \, \Omega \ = \ - \, 6 \, * \T_2 \, - \, 2 \, * (X \wedge \Omega) \, 
+ \, 2 \, * \T_{12} 
\edm 
we can express the difference $d \, \Omega \, - \, * \T^{\mathrm{c}}$ ,
\bdm
d \, \Omega \, - \, * \T^{\mathrm{c}} \ = \ 4 \cdot (* \, \T_{12} \, - \, 
* \, \T_{2}) \, . 
\edm
Consequently, we obtain
\begin{prop} 
The characteristic torsion form $\T^{\mathrm{c}}$ of a
$\mathrm{G}_1$-manifold is coclosed, $\delta(\T^{\mathrm{c}}) = 0$, if and only
if
\bdm
d \, * \T_2 \ = \ d \, * \T_{12}  .
\edm
In particular, any almost hermitian manifold of pure type  $\mathcal{W}_1$,
of pure type  $\mathcal{W}_3$ or of pure type  $\mathcal{W}_4$ has a coclosed 
characteristic torsion form. 
\end{prop}
%
\subsection{The characteristic connection of a nearly K\"ahler manifold}\noindent
\vspace{3mm}

\noindent
Nearly K\"ahler manifolds in dimension six have certain special properties.
They are Einstein spaces of positive scalar curvature, the almost complex
structure is never integrable, the first Chern class vanishes and they admit
a spin structure (see \cite{Gray76}). Moreover, nearly K\"ahler manifolds in 
dimension six are exactly those Riemannian spaces admitting real Riemannian 
Killing spinors (see \cite{FriedrichGrunewald}, \cite{Grunewald}). From our 
point of view, one of the interesting properties of nearly K\"ahler 
$6$-manifolds is the  $\nabla^{\mathrm{c}}$-parallelism of their torsion 
form $\T^{\mathrm{c}}$. This is a consequence of certain 
curvature identities already proved
by Takamatsu (see \cite{Takamatsu}), Matsumoto (see \cite{Matsumoto})
and Gray (see \cite{Gray76}). In Kirichenko's paper 
(\cite{Kirichenko}), the $\nabla^{\mathrm{c}}$-parallelism of 
$\T^{\mathrm{c}}$ appeared probably for the first time explicitly. 
We will outline a simple proof of this theorem. A nearly K\"ahler structure 
is characterized by the conditions
\bdm
Z \ = \ Z^* \, \oplus \, \Gamma_2 \, , \quad \Gamma_2 \ = \ \Theta(\T_2) \, ,
\quad \T_2 \in \Lambda_2^3  .
\edm
The derivative of the K\"ahler form and the characteristic torsion are
given by the formulas
\bdm
d * \Omega \ = \ 0 \, , \quad d \, \Omega \ = \ - \, 6 \, * \T_2 \, , \quad \T^{\mathrm{c}} \ = \ 
- \, 2 \, \T_2 \, , \quad \T^{\mathrm{c}}(X , Y) \ = \ - \,  
\J(\nabla^{\mathrm{LC}}_X\J)(Y) .
\edm
A nearly K\"ahler $6$-manifold is of constant type in the sense of
Gray (see \cite{Gray70}), i.e.,  
\bdm
\big|\big| \J(\nabla^{\mathrm{LC}}_X\J)(Y)\big|\big|^2 \ = \frac{\Scal}{30} \, 
\big\{\big|\big| X \big|\big|^2 \big|\big|Y\big|\big|^2 \, - \, g^2(X,Y) \, 
- \, g^2(X , \J(Y)) \big\} \, .
\edm
In particular, the length of the characteristic torsion coincides with
the scalar curvature,
\bdm
\big|\big| \T^{\mathrm{c}} \big| \big|^2 \ = \ \frac{2}{15} \, \Scal .
\edm
Since a nearly K\"ahler $6$-manifold is Einstein, the length of the characteristic
torsion is hence constant. It is a remarkable fact that this property  
of the characteristic connection implies alone that it is parallel.
\begin{thm} \label{nearlytorsion} 
The torsion of the characteristic connection of a nearly
K\"ahler $6$-manifold is parallel,
\bdm
\nabla^{\mathrm{c}} \T^{\mathrm{c}} \ = \ 0 \, .
\edm
The characteristic connection of a $6$-dimensional nearly K\"ahler 
non-K\"ahler manifold is a $\SU(3)$-connection.
\end{thm}
\begin{proof}
First of all, we remark that, for a $3$-form $\T_2 \in \Lambda_2^3$, we have
$\Gamma_2(X) =  X \haken \T_2$.
This implies that the characteristic connection coincides
with the connection $Z^*$ in the decomposition $Z = Z^* \oplus \Gamma_2$
of the Levi-Civita connection. The characteristic connection induces a metric
covariant derivative in the $2$-dimensional bundle $\Lambda_2^3$. Since $\T_2$ has constant length, there exists a $1$-form $A$ such that
\bdm
\nabla^{\mathrm{c}}_X \T_2 \ = \ A(X) \cdot (* \T_2) \, .
\edm
The codifferentials $\delta (\T) = \delta^{\nabla}(\T)$
of the torsion form of a metric connection coincide (see \cite{AgFr1}).
Therefore we obtain
\bdm
0 \, = \, \delta (* \, d \, \Omega) \, = \, 
\delta (6 \, \T_2) \, = \, - \, 3 \, 
\delta (\T^{\mathrm{c}}) \, = \, 
- \, 3 \, \delta^{\nabla^\mathrm{c}}(\T^{\mathrm{c}}) \, = \, 
6 \, \delta^{\nabla^\mathrm{c}}(\T_2) \, = \, 6 \, A \haken (* \T_2) \, .
\edm
The algebraic type $\T_2 \in \Lambda_2^3$ implies that $A = 0$ vanishes.
\end{proof}
\noindent
We remark that the K\"ahler form of a nearly K\"ahler $6$-manifold
is an eigenform of the Hodge-Laplace operator. Indeed, 
we can write the
equation $\nabla^{\mathrm{c}} \T^{\mathrm{c}} = 0$ equivalently as
\bdm
\nabla^{\mathrm{LC}}_X (* d \, \Omega) \, + \, \frac{\Scal}{10} \, 
\J(X) \wedge \Omega \ = \ 0 \, .
\edm
The latter formula immediately implies that
$5 \cdot \Delta \, \Omega \, = \, 2 \cdot \Scal \cdot \ \Omega$ .
\begin{NB}
The equation 
(see \cite[page 146-149]{YanoKon} or \cite[Theorem 5.2]{Gray76})
\bdm
\sum_{i=1}^6 g \big( \mathrm{R}(X , Y) e_i \, , \,  \J(e_i)\big) \ 
= \ - \, \frac{\Scal}{15} \, \Omega(X \, , \, Y)
\edm
is equivalent to the fact that the characteristic connection is an 
$\SU(3)$-connection. Indeed, the structure equation reads as
\bdm
\Omega^{\mathrm{LC}} \ = \ \Omega^{Z^*} \, + \, d \Gamma_2 \, + \, 
\big[Z^* , \, \Gamma_2 \big] \, + \,  \frac{1}{2} \big[ \Gamma_2 \, , 
\, \Gamma_2 \big] \, .
\edm
We project onto
the central element of the Lie algebra $\un(3)$. Since $d \Gamma_2$
and $[Z^* , \, \Gamma_2]$ have values in the subspace $\m^6$, we obtain
\bdm
\mathrm{pr}\big(\Omega^{\mathrm{LC}}\big) \ = \ 
\mathrm{pr} \big(\Omega^{Z^*}\big) \, + \,  \mathrm{pr} \big(\frac{1}{2} 
\big[ \Gamma_2 \, , \, \Gamma_2 \big]\big) \, .
\edm
The curvature identity of a nearly K\"ahler manifold mentioned above
as well as $ 15 \, || \T^{\mathrm{c}}||^2 = 2 \, \Scal$ yield that
\bdm
\mathrm{pr}\big(\Omega^{\mathrm{LC}}\big) \ = \ \mathrm{pr} \big(\frac{1}{2} 
\big[ \Gamma_2 \, , \, \Gamma_2 \big]\big) \ ,
\edm
i.e., the characteristic connection is a
$\SU(3)$-connection. 
\end{NB} 
\begin{NB}
The complete nearly K\"ahler manifolds with 
characteristic holonomy group contained in $\U(2) \times \U(1) \subset \U(3)$
have been
classified in \cite{BelgunMoroianu}. There are only two spaces of that type,
namely the projective space $\C \mathbb{P}^3$ and the flag manifold 
$\mathrm{F}(1,2)$ equipped with their homogeneous (non-K\"ahler) nearly 
K\"ahler structure. However, there is another interesting case. 
The $3$-dimensional complex irreducible representation of the group $\SU(2)/
\{\pm 1 \}$ 
is {\it reducible} as a real representation (see the discussion after
Theorem \ref{hologroup}). It is realized as the
characteristic holonomy by a left-invariant nearly K\"ahler structure on the 
Lie group $\mathrm{S}^3 \times \mathrm{S}^3$.
\end{NB}
\begin{NB} Homogeneous nearly K\"ahler manifolds have been
classified in \cite{Butruille}. The geometry of these examples has been
described in detail in \cite{BFGK}.
\end{NB} 
%
%
\subsection{$\mathrm{G}_1$-manifolds with parallel torsion and 
non-vanishing divergence}\noindent
\vspace{3mm}

\noindent
The aim of the next two sections is to study the structure of almost 
hermitian manifolds with a $\nabla^{\mathrm{c}}$-parallel
characteristic torsion $\T^{\mathrm{c}}$. We already know that any
nearly K\"ahler manifold has this property. Moreover, naturally reductive,
almost hermitian manifolds are automatically of type
$\mathrm{G}_1$ and their torsion form is $\nabla^{\mathrm{c}}$-parallel, too. 
Indeed, the canonical connection $\nabla^{\mathrm{can}}$  of a naturally
reductive space has totally skew-symmetric torsion and
preserves the almost hermitian structure. Since these two properties single out
the characteristic connection of the almost hermitian structure, 
we conclude that $\nabla^{\mathrm{c}}$ and $\nabla^{\mathrm{can}}$ coincide. 
But the canonical connection of any naturally reductive space has parallel 
torsion. This series of examples includes compact
Lie groups equipped with a left invariant almost hermitian structure.  
On the other side, left invariant almost complex structures
on nilmanifolds in dimension 
six have been discussed in detail in the papers 
\cite{AGS} and
\cite{FinoPortonSalamon}. Here the torsion form is, in general, not parallel.\\

\noindent 
In this section
we study $\mathrm{G}_1$-manifolds with a $\nabla^{\mathrm{c}}$-parallel
torsion form and non vanishing divergence of the K\"ahler form. 
The intrinsic torsion of a $\mathrm{G}_1$-manifold is given by two $3$-forms $\T_2 \, , \, \T_{12}$ and a $1$-form
$X$. The equations are
\bdm
d \, \Omega \ = \ - \, 6 \, * \T_2 \, - \, 2 \, * (X \wedge \Omega) \, + \, 
2 \, * \T_{12} \, , \quad \T^{\mathrm{c}} \ = \ - \, 2 \, \T_2 \, - \, 
2 \, (X \wedge \Omega) \, - \, 2 \, \T_{12} . 
\edm 
Since $\nabla^{\mathrm{c}} \, \T^{\mathrm{c}} = 0$ implies 
$\delta \T^{\mathrm{c}}  = \delta^{\nabla^{\mathrm{c}}}\T^{\mathrm{c}} = 0$,
we obtain the necessary conditions
\bdm
d \, \big(*\T_{12} \, - \, * \T_2 \big) \ = \ 0 \, , \quad
d \, \big( * (X \wedge \Omega) \, + \, 2 \, * \T_2 \big) \ = \ 0  .
\edm
The characteristic connection preserves the splitting $\Lambda^3 = 
\Lambda^3_2 \oplus \Lambda^3_6 \oplus \Lambda^3_{12}$. Therefore, the condition
$\nabla^{\mathrm{c}}\, \T^{\mathrm{c}} = 0$ is equivalent to
\bdm
\nabla^{\mathrm{c}} \, X \ = \ 0 \, , \quad
\nabla^{\mathrm{c}} \, \T_2 \ = \ 0 \, , \quad 
\nabla^{\mathrm{c}} \, \T_{12} \ = \ 0 . 
\edm
The forms $\mathrm{S}_1 := \T_{12} - \T_2$ and $\mathrm{S}_2 := X \wedge
\Omega + 2 \, \T_{12}$ are $\nabla^{\mathrm{c}}$-parallel and divergence free,
\bdm
\nabla^{\mathrm{c}} \, \mathrm{S}_1 \ = \ 0 \ = \ \nabla^{\mathrm{c}} \, 
\mathrm{S}_2 \, ,
\quad \delta \, \mathrm{S}_1 \ = \ 0 \ = \ 
\delta \, \mathrm{S}_2 .
\edm
Using the formula in \cite[Proposition 5.1]{AgFr1} we conclude that 
$(\alpha = 1,2)$
\bdm
\sum_{i,j=1}^6 \big(e_i \haken e_j \haken \T^{\mathrm{c}} \big) \wedge
\big( e_i \haken e_j \haken \mathrm{S}_{\alpha} \big) \ = \ 0  .
\edm
The latter equation couples the $3$-forms $\T_2$ and $\T_{12}$ via the
form $X$. Lenghtly, but elementary computations allow us to express this
link directly.
\begin{prop} \label{link}  
Let $M^6$ be a $\mathrm{G}_1$-manifold with parallel 
characteristic torsion, $\nabla^{\mathrm{c}} \, \T^{\mathrm{c}} = 0$. Then,
for any vectors $Y,Z$, the following equations are satisfied:
\begin{eqnarray*} 
\T_{12}(X \, , \, \J X \, , \, Y ) &=& 0 \, , \\
\T_{12}(X \, , \, Y \, , \, Z ) &=& \T_{12}( X \, , \, \J Y \, , \, \J Z ) 
\, + \, 2 \, \T_2(X \, , \, Y \, , \, Z) \, , \\
\T_{12}(\J X \, , \, Y \, , \, Z ) &=& \T_{12}(\J X \, , \, \J Y \, , \, 
\J Z ) \, - \, 2 \, \T_2(\J X \, , \, Y \, , \, Z) \ .
\end{eqnarray*}
\end{prop}
\noindent
The vector fields $X$ and $\J X$ are $\nabla^{\mathrm{c}}$-parallel and we 
compute their commutator,
\bdm
\big[X \, , \, \J X \big] \ = \ - \, \T^{\mathrm{c}}(X\,,\,\J X \,,\,-)\,.
\edm
For algebraic reasons, we have  $\T_2(X\,,\,\J X\,,\, -) = 0$ and
the formula simplifies
\bdm
\big[X \, , \, \J X \big] \ = \ 2 \, \T_{12}(X \, , \, \J X \, , \, -)\,.
\edm
The first equation of Proposition \ref{link} yields now the proof of the 
following
\begin{cor} 
Let $M^6$ be a $\mathrm{G}_1$-manifold with $\nabla^{\mathrm{c}}$-parallel 
characteristic torsion. Then $X$ and $\J X$ are commuting Killing
vector fields of constant length,
\bdm
\big[X \, , \, \J X \big] \ = \ 0 \, .
\edm
In case that the vector field $X \neq 0$ is non trivial, the leaves of the
integrable distribution $\{X \, , \, \J X \}$ are two-dimensional flat
and totally geodesic submanifolds. They are the orbits of an 
isometric $\R^2$-action.
\end{cor} 
\noindent
From now on we assume that the vector field $X \neq 0$ is non trivial. Then 
the tangent bundle splits into the integrable distribution 
$T^{\mathrm{v}} = \mathrm{Lin}(X \, , \, \J X)$ and its orthogonal complement 
$T^{\mathrm{h}}$. We decompose the $3$-forms $\T_2$ and $\T_{12}$ into
\bdm
\T_2 \ = \ X \wedge \Omega_1 \, + \, \J X \wedge \Omega_2 \, , \quad
\T_{12} \ = \ X \wedge \Omega_3 \, + \, \J X \wedge \Omega_4 .
\edm
$\Omega_1, \, \ldots , \, \Omega_4 \in \Lambda^2(T^{\mathrm{h}})$ are
horizontal $2$-forms. Remark that, for purely algebraic reasons, these forms
are orthogonal to the horizontal K\"ahler form $e_3 \wedge e_4 + e_5 \wedge 
e_6$. Proposition \ref{link} can be reformulated as
\bdm
\Omega_3 \ = \ \J (\Omega_3) \, + \, 2 \, \Omega_1 \, , \quad 
\Omega_4 \ = \ \J (\Omega_4) \, + \, 2 \, \Omega_2 .
\edm
and all these forms are $\nabla^{\mathrm{c}}$-parallel.
The next Proposition summarizes the result of a straightforward calculation.
\begin{prop} \label{Ableitung}  
The Lie derivative of the K\"ahler form and the differentials
of the forms $X$ and $\J X$ are given by
\begin{eqnarray*}
d X &=& ||X||^2 \, \big( \J (\Omega_3) \, - \, 3 \, \Omega_3 \, - \, 
2 \, \Omega \big) \, + \, 2 \, X \wedge \J X \, , \\ 
d \, \J X &=& ||X||^2 \, \big( \J (\Omega_4)\,-\, 3 \,\Omega_4 \big)\,,\quad
\mathcal{L}_X \Omega \ = \ 8 \, ||X||^2 \, \Omega_2 .
\end{eqnarray*} 
\end{prop}
\noindent
A direct consequence of these formulas is
\begin{thm} 
Let $(M^6, g, \J)$ be an almost hermitian $6$-manifold of type
$\mathcal{W}_1 \oplus \mathcal{W}_4$. If the torsion of its characteristic
connection is parallel, $\nabla^{\mathrm{c}} \, \T^{\mathrm{c}} = 0$, then
$M^6$ is either of pure type $\mathcal{W}_1$ or of
pure type $\mathcal{W}_4$.
\end{thm}
\begin{NB} 
The characteristic torsion of a nearly K\"ahler manifold 
is $\nabla^{\mathrm{c}}$-parallel. On the other hand, 
suppose that $M^6$ is of pure type
$\mathcal{W}_4$ and $X$ is $\nabla^{\mathrm{c}}$-parallel. Then we obtain
\bdm
0 \ = \ \nabla^{\mathrm{c}}_Z \, \J X \ = \ 
\nabla^{\mathrm{LC}}_Z \, \J X \, , \quad d X \ = \ 2 \, \big( -  ||X||^2 
\, \Omega \, + \, X \wedge \J X \big) \, .  
\edm
The vector field $\J X$ is $\nabla^{\mathrm{LC}}$-parallel, i.e., the manifold 
is a generalized Hopf manifold. Up to a 
scaling of the length of $X$,  the manifold is 
locally isometric to a product $M^6 = N^5 \times \R^1$ of 
the line $\R^1$ by a $5$-dimensional Sasakian manifold $N^5$. Conversely, 
any product of a Sasakian manifold by $\R^1$ is an almost hermitian manifold of
type $\mathcal{W}_4$ with parallel torsion. These  $\mathcal{W}_4$-manifolds
have been studied by Vaisman, see \cite{Vaisman}.
\end{NB}
\noindent
We consider now hermitian manifolds. The complex structure is integrable and
the forms $\Omega_1 $ and $\Omega_2 $ vanish. The forms 
$\Omega_3, \, \Omega_4 \in \Lambda^2_-(T^{\mathrm{v}})$ are anti-selfdual with
respect to the $4$-dimensional horizontal Hodge operator and the formulas of 
Proposition \ref{Ableitung} simplify,
\bdm
d X \ = \ - \, 2 \ ||X||^2 \, \Omega_3 \, - \, 
2 \, ||X||^2 \, \Omega  \, + \, 2 \, X \wedge \J X \, , \quad 
d \, \J X \ = \ - \, 2 \, ||X||^2 \, \Omega_4 .
\edm
Furthermore, we obtain immediately that
\bdm
\mathcal{L}_X \Omega \ = \ 0 \, , \quad 
\mathcal{L}_{\J X} \Omega \ = \ 0 \, , \quad 
d \big( ||X||^2 \, \Omega \, - \, 
X \wedge \J X \big) \ = \ 0  .
\edm
The forms $\Omega_3$ and $\Omega_4$ are closed. Suppose that the 
Killing vector fields $X$ and $\J X$ induce a regular group action, i.e., the
orbit space $\hat{X}^4$ is smooth. Then $\hat{X}^4$ admits a Riemannian
metric $\hat{g}$ and a complex structure $\hat{\J}$ with K\"ahler form
\bdm
\hat{\Omega} \ = \ \Omega \ - \ \frac{1}{||X||^2} \, X \wedge \J X \ .
\edm
In particular, $\hat{X}^4$ is a $4$-dimensional K\"ahler manifold.
The forms $\Omega_3$ and $\Omega_4$ project onto $\hat{X}^4$. A direct
computation yields that the $\nabla^{\mathrm{c}}$-parallelism of these
forms on $M^6$ can be reformulated as the condition that their projections 
$\hat{\Omega}_3$ and $\hat{\Omega}_4$ are anti-selfdual and parallel forms
on the K\"ahler manifold $\hat{X}^4$, 
\bdm
\hat{\nabla} \, \hat{\Omega}_3 \ = \ 0 \, , \quad 
\hat{\nabla} \, \hat{\Omega}_4 \ = \ 0 \, , \quad
* \, \hat{\Omega}_3 \ = \ - \, \hat{\Omega}_3 \, , \quad
* \, \hat{\Omega}_4 \ = \ - \, \hat{\Omega}_4  .
\edm
The structure group of the principal fiber bundle $M^6 \rightarrow \hat{X}^4$
is $2$-dimensional and abelian. Up to a scaling of the length, 
the pair $\{X \, , \, \J X \}$ is a connection. Its curvature is the pair of
$2$-forms 
$(\mathrm{Curl}_1 \, , \mathrm{Curl}_2)$ on $\hat{X}^4$ given by the
differentials of $X$ and $\J X$, i.e.
\bdm
\mathrm{Curl}_1 \ = \ - \, 2 \, \hat{\Omega}_3 \ - \ 2 \, \hat{\Omega} \, , 
\quad \mathrm{Curl}_2 \ = \ - \, 2 \, \hat{\Omega}_4 .
\edm
Vice versa, we can reconstruct the whole $6$-dimensional structure out of
the $4$-dimensional K\"ahler manifold $(\hat{X}^4 \,,\,\hat{g}\,,\,\hat{\J})$ 
and the two parallel forms $\hat{\Omega}_3 \, , \, \hat{\Omega}_4 
\in \Lambda^2_-(\hat{X}^4)$. In the compact case we need that 
$2 \, \hat{\Omega}_4$ and $2 \, \hat{\Omega} \, + \, 2 \, \hat{\Omega}_3$ are 
curvature forms of some $\U(1)$-bundle, i.e. , 
\bdm 
2 \, \hat{\Omega}_4 \, , \   2 \, \hat{\Omega} \, + \, 2 \, \hat{\Omega}_3
 \in \mathrm{H}^2(\hat{X}^4 \, ; \Z) \ .
\edm
We summarize the result for compact hermitian spaces $M^6$.
\begin{thm}
The compact regular hermitian manifolds $(M^6 \, , \, g \, , \, \J)$ 
with non vanishing divergence  $- \, 4 \, X = \delta \, \Omega \neq 0$ of the 
K\"ahler form and $\nabla^{\mathrm{c}}$-parallel characteristic torsion  
$\T^{\mathrm{c}}$ correspond to triples $(\hat{X}^4 \, , \, \hat{\Omega}_3 
\, , \, \hat{\Omega}_4)$ consisting of a compact $4$-dimensional K\"ahler 
manifold $\hat{X}^4$ and two parallel anti-selfdual forms $\hat{\Omega}_3$ ,
$\hat{\Omega}_4$ such that
\bdm 
2 \, \hat{\Omega}_4 \, , \   2 \, \hat{\Omega} \, + \, 2 \, \hat{\Omega}_3
 \in \mathrm{H}^2(\hat{X}^4 \, ; \Z) \ .
\edm   
\end{thm} 
\noindent
It is easy to describe the possible K\"ahler manifolds $\hat{X}^4$. First
of all, a parallel anti-selfdual $2$-form gives rise to a parallel complex
structure of opposite orientation. Then a compact $4$-dimensional space with
two independent parallel complex structures with equal orientation is
hyperk\"ahler. The existence of opposite parallel complex structures restricts
it to be a torus (see \cite{Hitchin74}). Since toric bundles over tori are
always $2$-step nilmanifolds, the $6$-dimensional manifold is at least
diffeomorphic to a locally homogeneous space. When the forms $\hat{\Omega}_3$
and $\hat{\Omega}_4$  are linearly depent, $M^6$ is actually a product
$\mathrm{S}^1 \times N^5$ and $N^5$ is a $\mathrm{S}^1$-bundle over
a $4$-manifold which is covered by a  product of two 
surfaces\footnote{The authors would like to thank the referee for a hint 
completing this classification.}. 
\begin{NB} 
Hermitian structures with a $\SU(3)$-holonomy of the characteristic
connection have been constructed recently on certain toric bundles,
see \cite{GrantcharovPonn}. The condition $\nabla^{\mathrm{c}}\T^{\mathrm{c}}
= 0$ is a stronger condition and, consequently, our family is much smaller.
\end{NB}
%
\subsection{$\mathcal{W}_3$-manifolds with parallel torsion}\noindent
\vspace{3mm}

\noindent
An interesting problem is the structure of $\mathcal{W}_3$-manifolds with
$\nabla^{\mathrm{c}}$-parallel torsion. The equations characterizing these
hermitian manifolds are (see \cite{FriedrichIvanov})
\bdm
d \, \Omega \ = \ 2 \, * \, \T_{12} \ = \ - \, * \, \T^{\mathrm{c}} \, , 
\quad  \delta \, \Omega \ = \ 0 \, , \quad \nabla^{\mathrm{c}} \T^{\mathrm{c}}
\ = \ 0 \, , \quad d\T^{\mathrm{c}} \ = \  2 \, \sigma_{\T^{\mathrm{c}}} \, ,
\edm
where the $4$-form $\sigma_{\T^{\mathrm{c}}}$ is defined by the formula
\bdm
\sigma_{\T^{\mathrm{c}}} \ := \ \frac{1}{2}\sum_{i=1}^6 (e_i \haken 
\T^{\mathrm{c}}) \wedge (e_i \haken \T^{\mathrm{c}}) \, .
\edm
We remark that in the class of hermitian  $\mathcal{W}_3$-manifolds an 
analogue of Theorem \ref{nearlytorsion} does not hold.
\begin{exa} 
Consider the $3$-dimensional complex Heisenberg group. There
exists a left-invariant metric with the following structure equations
\bdm
d e_1 \ = \ d  e_2 \ = \ d  e_3 \ = \ d e_4 \ = \ 0 \, , \quad
d e_5 \ = \ e_{13} \, - \, e_{24} \, , \quad
d e_6 \ = \ e_{14} \, + \, e_{23} .
\edm 
The differential of the K\"ahler form is given by
\bdm
d \, \Omega \ = \ e_{136} \, - \, e_{246} \, - \, e_{145} \, - \, e_{235} .
\edm
Consequently, the hermitian structure is of pure type
$\mathcal{W}_3$ and its torsion is given by $\T^{\mathrm{c}} = 
e_{245} \, - \, e_{135} \, - \, e_{236} \, - \, e_{146}$. We compute
the derivative $d \T^{\mathrm{c}}$ and the $4$-form
$\sigma_{\T^{\mathrm{c}}}$,
\bdm
d \, \T^{\mathrm{c}} \ = \ - \, 4 \, e_{1234} \, , \quad
\sigma_{\T^{\mathrm{c}}} \ = \ 2 \, e_{1234} \, - \, (e_{12} \, + \, e_{34}) 
\wedge e_{56} . 
\edm
Since $d \T^{\mathrm{c}} \neq 2 \, \sigma_{\T^{\mathrm{c}}}$, 
the torsion form of the Heisenberg group is not parallel.
\end{exa} 
\noindent
The $\U(3)$-orbit type of the parallel torsion form 
$\T^{\mathrm{c}} \in \Lambda^3_{12}$ is constant. There are only two types of 
$3$-forms in $ \Lambda^3_{12}$ with a non abelian isotropy group.
\begin{thm} \label{hologroup}  
Let $\T \in \Lambda^3_{12}$ be a $3$-form and denote by
$\mathrm{G}_{\T} \subset \U(3)$ the connected component of its
isotropy group. If the dimension of $\mathrm{G}_{\T}$ is at least three, 
then one of
the following two cases occurs.
\begin{enumerate}
\item
The group $\mathrm{G}_{\T}$ is isomorphic to $\U(2)$ and the embedding
into $\U(3)$ is given by the homomorphism
\bdm
\mathrm{G}_{\T} \ = \ 
\Big \{ \left[\ba{cc} g & 0 \\ 
0 & \mathrm{det}(g) \ea\right]  , \ g \in \U(2) \Big \} \, .
\edm
Up to a complex factor, there exists one orbit of that type represented by 
the $3$-form
\bdm
\T \ = \  (e_{135} \, - \, e_{245} \, + \, e_{236} \, + \, e_{146})\, .
\edm
\item
The group $\mathrm{G}_{\T}$ is isomorphic to $\SU(2)/\{\pm 1\} = \SO(3)$ 
and the embedding
into $\U(3)$ is the unique $3$-dimensional irreducible complex
representation of $\SU(2)$. Up to a complex factor, there
exists one orbit of that type represented by the $3$-form
\bdm
\T \ = \ 2 \, (e_{123} \, - \, e_{356}) \, - \, (e_{246} \, + \, e_{136}) 
\, + \, (e_{145} \, - \, e_{235}) \, .
\edm
\end{enumerate}  
\end{thm}
\begin{proof}
We use the explicit equations defining the Lie algebra 
$\g_{\T} \subset \un(3)$ of the isotropy group  $\mathrm{G}_{\T}$. 
The $3$-form $\T$ depends on $12$ real parameters,
\begin{eqnarray*}
\T &=& A_1 \, (e_{123} \, - \, e_{356}) \, + \, A_2 \, (e_{124} \, - 
\, e_{456}) \, + \, A_3 \, (e_{125} \, - \, e_{345}) \, + \, 
A_4 \, (e_{126} \, - \, e_{346}) \, + \\
&& A_5 \, (e_{134} \, - \, e_{156}) \, + \, 
A_6 \, (e_{234} \, - \, e_{256}) \, + \,
A_7 \, (e_{135} \, + \, e_{245}) \, + \, A_8 \, 
(e_{246} \, + \, e_{136}) \, + \\ 
&& A_9 \, (e_{145} \, - \, e_{235}) \, + \, A_{10} \, (e_{236} \, - \, 
e_{146}) \, + \, A_{11} \, (e_{135} \, - \, e_{245} \, + \, e_{236} \, + \, 
e_{146}) \, + \\
&& A_{12} \, (e_{246} \, - \, e_{136} \, + \, e_{145} \, + \, e_{235}) \, . 
\end{eqnarray*}
An arbitrary $2$-form in $\un(3)$ depends on $9$ real parameters,
\begin{eqnarray*}
\omega &=& \omega_{12} \, e_{12} \, + \, 
\omega_{13} \, (e_{13} \, + \, e_{24}) \, + \, 
\omega_{14} \, (e_{14} \, - \, e_{23}) \, + \,
\omega_{15} \, (e_{15} \, + \, e_{26}) \, + \\
&& \omega_{16} \, (e_{16} \, - \, e_{25}) \, + \,
\omega_{34} \, e_{34} \, + \, 
\omega_{35} \, (e_{35} \, + \, e_{46}) \, + \,
\omega_{36} \, (e_{36} \, - \, e_{45}) \, + \,  
\omega_{56} \, e_{56} .
\end{eqnarray*}
The condition $\vrho_*(\omega) \, \T \, = \, 0$ is a linear system of $12$ 
equations with respect to $9$ variables $\omega_{ij}$ given by the 
following $(12 \times 9)$--matrix $\mathcal{A}_{\T}$,
\bdm
\left[\ba{ccccccccc} \kl{2 A_{12}} & \kl{0} & \kl{0} & \kl{2 A_2} & 
\kl{- 2 A_1} & \kl{2 A_{12}} & \kl{- 2 A_6} &
\kl{2 A_5} & \kl{- 2 A_{12}} \\ 
\kl{-2A_{11}} & \kl{0} & \kl{0} & \kl{-2A_1} & \kl{-2A_2} & 
\kl{-2A_{11}} & \kl{2A_5} & \kl{2A_6} & \kl{2A_{11}} \\  
\kl{A_{6}} & \kl{A_1} & \kl{A_2} & \kl{-A_3} & \kl{-A_4} & \kl{0} 
& \kl{D - 2 A_{12}} & \kl{- B - 2A_{11}} & \kl{0} \\ 
\kl{-A_5} & \kl{A_2} & \kl{-A_1} & \kl{-A_4} & \kl{A_3} & \kl{0} 
& \kl{- B + 2A_{11}}  & \kl{- D - 2 A_{12}} & \kl{0}\\
\kl{0} & \kl{A_5} & \kl{-A_6} &  \kl{C - 2 A_{12}} & \kl{- A - 2A_{11}} & 
\kl{-A_2} & \kl{-A_3} & \kl{-A_4} & \kl{0} \\
\kl{0} & \kl{-A_6} & \kl{-A_5} &  \kl{A - 2 A_{11}} & 
\kl{C + 2A_{12}} & \kl{-A_1} & \kl{A_4} & \kl{-A_3} & \kl{0} \\
\kl{- D} & \kl{2A_4} & \kl{-2A_3} & \kl{0} & \kl{0} & \kl{D} & 
\kl{-2 A_6} & \kl{-2 A_5} & \kl{- D} \\
\kl{B} & \kl{-2A_3} & \kl{-2A_4} & \kl{0} & \kl{0} & \kl{- B} 
& \kl{2 A_5} & \kl{-2 A_6} & \kl{B} \\
\kl{0} & \kl{2A_{10}} & \kl{2 A_8} & \kl{-A_6} & \kl{-A_5} 
& \kl{0} & \kl{A_2} & \kl{A_1} & \kl{A_3} \\
\kl{0} & \kl{-2A_{9}} & \kl{2 A_7} & \kl{-A_5} & \kl{A_6} 
& \kl{0} & \kl{A_1} & \kl{-A_2} & \kl{-A_4} \\
\kl{C} & \kl{- 2A_4} & \kl{- 2A_3} & \kl{2A_2} & \kl{2A_1} 
& \kl{- C} & \kl{0} & \kl{0} & \kl{- C}\\
\kl{A} & \kl{- 2A_3} & \kl{2A_4} & \kl{2A_1} & \kl{-2A_2} 
& \kl{- A} & \kl{0} & \kl{0} & \kl{- A} \ea\right] \, . 
\edm
We introduced the notion $A := A_7 + A_{10} , \, 
B := A_7 - A_{10} , \, C := A_8 + A_9 , \, D := A_8 - A_9$.
If $\T \neq 0$,  the rank of this matrix is at least three. Therefore,
the dimension of the Lie algebra is bounded by $\mathrm{dim}(\g_{\T}) \leq
6$. Since $\T \in \Lambda^3_{12}$, the central element $\Omega \in \un(3)$ 
does not belong to $\g_{\T}$. An elementary discussion concerning
subgroups of $\U(3)$ yields the result that the group $\mathrm{G}_{\T}$
is conjugated to $\SO(3)$ or $\U(2)$ and realized in the way as the
theorem states. On the other side, given one of these 
two Lie algebras $\g_{\T}$, the matrix $\mathcal{A}_{\T}$ computes the 
corresponding $3$-form $\T$ up to a complex factor.
\end{proof}
\noindent
First, we study the case of $\mathrm{G}_{\T^{\mathrm{c}}} = \U(2)$. Then 
the $2$-forms $e_{12} + e_{34}$ and 
$e_{56}$ are globally defined and $\nabla^{\mathrm{c}}$-parallel,
\bdm
\nabla^{\mathrm{c}} (e_{12} \, + \, e_{34}) \ = \ 0 \ = \ \nabla^{\mathrm{c}} 
(e_{56}) \, .
\edm 
Using \cite[Proposition 5.2]{AgFr1}, we compute the exterior derivative 
$d(e_{56})$,
\bdm
d(e_{56}) \ = \ \sum_{i=1}^6 \big( e_i \haken e_{56} \big) \wedge \big(
e_i \haken \T^{\mathrm{c}} \big) \ = \ * \, \T^{\mathrm{c}} .
\edm  
Moreover, $d \Omega = - \, *  \T^{\mathrm{c}}$ implies a formula for the
derivative of the second invariant $2$-form, 
$d(e_{12} + e_{34}) = - \, 2 \, * \T^{\mathrm{c}}$. Let us introduce a new 
almost complex structure $\hat{\J}$ by the condition
\bdm
\hat{\Omega} \ = \ - \, (e_{12} \, + \, e_{34}) \, + \, e_{56} .
\edm
Then we have
\bdm
\nabla^{\mathrm{c}} \hat{\Omega} \ = \ 0 \, , \quad 
d \, \hat{\Omega} \ = \ 3 \, * \, \T^{\mathrm{c}} \ ,
\edm 
i.e., the manifold $(M^6 , g ,  \hat{\J})$ is nearly K\"ahler, 
$\nabla^{\mathrm{c}}$ is its characteristic connection, and the holonomy 
$\mathrm{Hol}(\nabla^{\mathrm{c}}) = \U(2) = \mathrm{G}_{\T}$ is not the whole
group $\SU(3)$. In the compact case these nearly K\"ahler manifolds have been
classified in \cite{BelgunMoroianu}. There are only two of them, namely
the twistor spaces of the $4$-dimensional sphere or of the complex projective
plain equipped with their canonical non-integrable almost complex structure and
their canonical non-K\"ahler Einstein metric. Replacing again
the almost complex structure $\hat{\J}$ by $\J$, we obtain  
a complete classification of all $\mathcal{W}_3$-manifolds with parallel 
characteristic torsion of type $\mathrm{G}_{\T^{\mathrm{c}}} = \U(2)$.
\begin{thm} 
Let $(M^6 , g , \J)$ be a complete hermitian manifold of type $\mathcal{W}_3$ 
such that
\bdm
\nabla^{\mathrm{c}} \T^{\mathrm{c}} \ = \ 0 \, , \quad 
\mathrm{G}_{\T^{\mathrm{c}}} \ = \ \U(2) \, .
\edm
Then $M^6$ is the twistor space of a $4$-dimensional, compact selfdual
Einstein manifold with positive scalar curvature. The complex structure
$\J$ is the standard one of the twistor space and 
the metric $g$ is the unique non-K\"ahler Einstein metric in the
canonical $1$-parameter family of metrics of the twistor space. 
\end{thm}
\begin{NB}
The latter Theorem holds locally and in higher dimensions too, 
see \cite{Alexandrov}. In dimension six, there are only two compact 
K\"ahlerian twistor spaces, namely the projective space $\CP^3$ and the flag
manifold $\mathrm{F}(1,2)$ (see \cite{FriedrichKurke} and \cite{Hitchin}). 
\end{NB}
\noindent
The second case $\mathrm{G}_{\T^{\mathrm{c}}} = \SU(2)/\{\pm 1\}\subset \U(3)$ 
corresponds to the $3$-dimensional {\it complex irreducible} representation.
The underlying real representation in $\C^3 = \R^6$ is {\it reducible}, i.e.,  
under the action of the group $\mathrm{G}_{\T^{\mathrm{c}}}$ the 
euclidian space $\R^6$ splits into two real and $3$-dimensional Lagrangian 
subspaces. The holonomy representation
is the sum of two faithful representations. The results
of \cite[Lemma 4.4 and Lemma 5.6]{CleytonSwann} yield that 
$M^6$ is a so-called Ambrose-Singer manifold,
i.e., the curvature $\mathrm{R}^{\mathrm{c}}$ of the characteristic
connection is $\nabla^{\mathrm{c}}$-parallel,
\bdm
\nabla^{\mathrm{c}} \T^{\mathrm{c}} \ = \ 0 \ , \quad 
\nabla^{\mathrm{c}} \mathrm{R}^{\mathrm{c}} \ = \ 0 \, .
\edm
Since the universal covering of  $\mathrm{G}_{\T^{\mathrm{c}}}$ is compact,
the Ambrose-Singer manifold is regular and locally isometric
to a homogeneous space $\mathrm{G}/\mathrm{G}_{\T^{\mathrm{c}}}$. 
The Lie algebra of the automorphism group $\mathrm{G}$ is the vector space
$\g := \g_{\T^{\mathrm{c}}} \oplus \R^6$ equipped with the bracket
(see \cite[Theorem 5.10]{CleytonSwann})
\bdm
\big[ A \, + \, X \, , \, B \, + \, Y \big] \ = \ 
\big( [A \, , \, B] \, - \, \mathrm{R}^{\mathrm{c}}(X , Y) \big) \, + \, 
\big( A \cdot Y \, - \, B \cdot X \, - \, \T^{\mathrm{c}}(X , Y) \big) \, .
\edm
In order to find the
automorphism group as well as the hermitian manifold, we consider 
the Lie subalgebra $\so(3) \subset \so(6)$. It is generated by
the following $2$-forms
\bdm
\omega_1 \, := \, \frac{1}{\sqrt{2}} \big( e_{12}  \, - \, e_{56} \big) \, , \,
\omega_2 \, := \, \frac{1}{2} \big( e_{13} \, + \, e_{24} \, + \, e_{35} \, +
\,  e_{46} \big) \, , \,
\omega_3 \, := \, \frac{1}{2} \big( e_{14} \, - \, e_{23} \, + \, e_{36} \, -
\,  e_{45} \big) \, ,
\edm
and the $\SO(3)$-invariant form $\T^{\mathrm{c}} \in
\Lambda^3(\R^6)$ is given by the formula
\bdm
\T^{\mathrm{c}} \ := \ 2 \big(e_{123} \, - \, e_{356}\big) \, - \, 
\big( e_{246} \, + \, e_{136} \, - \, e_{145} \, + \, 
e_{235} \big) \, .
\edm
The 
curvature tensor of the characteristic connection is an $\SO(3)$-invariant 
$2$-form with values in the Lie algebra  $\so(3)$. Since the
$\SO(3)$-representation $\Lambda^2(\R^6)$ splits into 
$3 \cdot \R^3 \oplus \R^1 \oplus \mathrm{S}_0^2(\R^3)$, the curvature tensor
depends a priori on three parameters. However, the first Bianchi identity yields that
$\mathrm{R}^{\mathrm{c}}$ is {\it unique}. 
\begin{lem}
The curvature of the characteristic connection is proportional to 
the orthogonal projection onto $\so(3)$,
\bdm
\mathrm{R}^{\mathrm{c}}  :  \Lambda^2(\R^6) \, = \, \so(6) \rightarrow 
\so(3) \, , \quad  \mathrm{R}^{\mathrm{c}}(X\, , \, Y) \ = \ 
4 \cdot \mathrm{pr}_{\so(3)}(X \wedge Y) \ .
\edm
\end{lem}
\noindent
We remark that the $3$-form $\T^{\mathrm{c}}$ 
satisfies the necessary condition 
in order to define an extension of the Lie algebra $\so(3)$, namely, 
the element of the Clifford algebra $\mathfrak{Cliff}(\R^6)$ 
\bdm
\big(\T^{\mathrm{c}}\big)^2 \, + \, 4 \cdot 
\big(\omega_1^2 \, + \, \omega_2^2 \, + \, \omega_3^2 \big) 
\edm
is a scalar (see \cite[chapter 10.4]{Sternberg}). It turns out that 
the automorphism group is isomorphic to the semi-simple Lie group 
$\mathrm{G} = \SL(2,\C) \times \SU(2)$. The hermitian manifold 
$M^6 = \mathrm{G}/\mathrm{G}_{\T^{\mathrm{c}}}$ is a left invariant hermitian
structure on $\SL(2,\C)$ represented as a naturally
reductive space by the help of the subgroup $\SU(2) \subset \SL(2,\C)$
(see \cite{Agri}, \cite{DAtriZiller}). Since the characteristic connection of 
the hermitian manifold is unique, it coincides with the canonical connection
of the naturally reductive space.
\begin{thm} 
Any hermitian $6$-manifold of type $\mathcal{W}_3$ and
\bdm
\nabla^{\mathrm{c}} \T^{\mathrm{c}} \ = \ 0 \, , \quad 
\mathrm{G}_{\T^{\mathrm{c}}} \ = \ \SU(2)/\{\pm 1\} 
\edm
is locally isomorphic to the left invariant hermitian structure
on the Lie group $\SL(2,\C)$.  
\end{thm}
\noindent
We briefly describe the hermitian structure under consideration.
Let us decompose the Lie algebra $\g = \slin(2,\C) \oplus \su(2)$,   
\bdm
\g \ = \ \big\{(A \, , \, B ) \in \M(2,\C) \oplus \M(2,\C) \, : \, 
\mathrm{tr}(A) \, = \, 0 \, , \ B \, + \, \overline{B}^t \, = \, 0 \, ,
\ \mathrm{tr}(B) \, = \, 0 \big\} \, ,
\edm
into the subalgebra $\h := \big\{(B \, , \, B ) \in \g \, : \, 
B \, + \, \overline{B}^t \, = \, 0 \, ,
\, \mathrm{tr}(B) \, = \, 0 \big\}$ and its complement,
\begin{eqnarray*}
\m \, := \, \big\{(A \, , \, B ) \in \g \, : \, A \, - \, 
\overline{A}^t \, = \, 0 \, , \,  \mathrm{tr}(A) \, = \, 0 \, , \, 
B \, + \, \overline{B}^t \, = \, 0 \, ,
\, \mathrm{tr}(B) \, = \, 0 \big\} \ .
\end{eqnarray*}
The decomposition is reductive, $[\h \, , \, \m ] \subset \m$. Moreover,
we introduce a complex structure $\J : \m \rightarrow \m$ as well as a scalar 
product $\langle \, , \, \rangle_{\m}$ by the formulas
\bdm
\J(A \, , \, B) \, := \, (i \cdot B \, , \, i \cdot A) \, , \quad
\langle (A \, , \, B) \, , \, (A_1 \, , \, B_1) \rangle_{\m} \, := 
\, \mathrm{tr}(A \cdot \overline{A}_1^t) \, + \, 
\mathrm{tr}(B \cdot \overline{B}_1^t) \, .
\edm
Both are $\h = \su(2)$-invariant. They define an almost hermitian structure 
on $M^6 = (\SL(2,\C) \times \SU(2))/\SU(2) = \SL(2,\C)$. It turns out that 
the almost complex structure is integrable
and the hermitian structure is of type $\mathcal{W}_3$ $( \delta \Omega = 0)$. 
Its characteristic torsion form coincides with the canonical torsion of the naturally reductive space. The manifold realizes the orbit type 
$\mathrm{G}_{\T^{\mathrm{c}}} = \SU(2)/\{\pm 1\}$. Finally, let us describe 
some geometric data. The Ricci tensor of the characteristic connection is 
proportional to the metric,
\bdm
\mathrm{Ric}^{\nabla^{\mathrm{c}}} \ = \ - \, \frac{1}{3} \cdot 
|| \T^{\mathrm{c}} ||^2 \cdot \mathrm{Id} \, .
\edm
The $3$-form $\T^{\mathrm{c}}$ acts on the spinor bundles
$S^{\pm}$ with a one-dimensional kernel and there exist
two $\nabla^{\mathrm{c}}$-parallel spinor fields $\Psi^{\pm}$,
\bdm
\nabla^{\mathrm{c}} \Psi^{\pm} \ = \ 0 \, , \quad \T^{\mathrm{c}} \cdot \Psi^{\pm} \ = \ 0 \, , \quad \nabla^{\mathrm{c}}\T^{\mathrm{c}} \ = \ 0 \, ,
\quad \delta(\T^{\mathrm{c}}) \ = \ 0 \, .
\edm 

\noindent
We study the case of $\mathrm{dim}\big(\mathrm{G}_{\T^{\mathrm{c}}}\big) 
\leq 2$
in a similar manner. Since $\mathrm{Hol}(\nabla^{\mathrm{c}}) \subset
\mathrm{G}_{\T^{\mathrm{c}}}$, we have the following possibilities:
\vspace{4mm}

\begin{center}
\begin{tabular}{|c|c|c|c|}
\hline 
$\mathrm{dim}\big(\mathrm{G}_{\T^{\mathrm{c}}}\big)$ & 0 , 1 , 2 
& 1 & 2 \\[1mm]
\hline
$\mathrm{dim}\big(\mathrm{Hol}(\nabla^{\mathrm{c}})\big)$ & 0 & 1 & 1 , 2 \\[1mm]
\hline
\end{tabular}
\end{center}
\vspace{4mm}

\noindent
If the holonomy group is
discrete, the characteristic connection is flat and
the manifold $M^6$ is a Lie group. Its Lie algebra is given by
$\g \, = \, \R^6 \, , \ 
\big[ X \, , \, Y \big] \, = \, - \, \T^{\mathrm{c}}(X , Y)$.
The Jacobi identity is equivalent to the condition that the square
$(\T^{\mathrm{c}})^2$ of the torsion form in the Clifford 
algebra $\mathfrak{Cliff}(\R^6)$ is a scalar (see \cite[Theorem 1.50]{Kostant} 
and 
\cite[chapter 10.4]{Sternberg}).  However, $3$-forms of type $\Lambda^3_{12}$
satisfying this condition do not exist. 
\begin{lem}\label{T^2}
Let $\T \in \Lambda^3_{12}$ be a $3$-form and such that its
square $\T^2 $ in $\mathfrak{Cliff}(\R^6)$ is a scalar. Then $\T = 0$.
\end{lem}
\begin{proof}
We parameterize a form $\T \in \Lambda^3_{12}$ by its coefficients 
$A_1, \cdots , A_{12}$
with respect to the introduced basis. The endomorphism $\T^2$
in  $\mathfrak{Cliff}(\R^6) \subset \mathfrak{Cliff}(\R^7) \rightarrow
\mathrm{End}(\Delta_7)$ acting in the $7$-dimensional real spin 
representation is an $(8 \times 8)$-matrix. We compute the numbers
on the diagonal :
\bdm
0 \, , \ 4 \cdot \big( A_1^2 + A_2^2 + A_5^2 + A_6^2 + A_{11}^2 + A_{12}^2
\big) \, , \ 
4 \cdot \big( A_3^2 + A_4^2 + A_5^2 + A_6^2 + (A_{7} \pm A_{10})^2 + 
(A_{8} \pm A_9)^2 \big) \, .
\edm 
Consequently,  $\T^2$ is a scalar if and only if $\T = 0$.
\end{proof}
\begin{thm} Let $(M^6, g, \J)$ be a complete hermitian manifold of type
$\mathcal{W}_3$ such that
\bdm
\nabla^{\mathrm{c}} \T^{\mathrm{c}} \ = \ 0 \, , \quad 
\mathrm{dim}\big(\mathrm{Hol}(\nabla^{\mathrm{c}})\big) \ = \ 0 \ . 
\edm
Then $M^6$ is a flat K\"ahler manifold, i.e. $\T^{\mathrm{c}} = 0$.
\end{thm}
\noindent
In the next step of our classification we will prove that 
$\mathrm{dim}\big(\mathrm{Hol}(\nabla^{\mathrm{c}})\big) = 1 =
\mathrm{dim}\big(\mathrm{G}_{\T^{\mathrm{c}}}\big)$ is impossible. 
The holonomy representation $ \mathrm{Hol}(\nabla^{\mathrm{c}}) =  
\mathrm{G}_{\T^{\mathrm{c}}} \rightarrow
\U(3)$ is given by three integers $k_1, k_2, k_3 \in \Z$ and the diagonal
matrices $ \varphi \rightarrow \mathrm{diag}( e^{i k_1 \varphi}, 
\ e^{i k_2 \varphi} , \ e^{i k_3 \varphi} )$.
If $k_2 , \, k_3 = 0$, the linear system $\rho_*(\omega) \T = 0$ has
a $4$-dimensional solution with respect to $\T$, namely
\bdm
A_5 \ = \ A_6 \ = \ A_7 \ = \ A_8 \ = \ A_9 \ = \ A_{10} \ = \ A_{11} 
\ = A_{12} \ = \ 0 \, .
\edm
However, a direct computation shows that for any of these $3$-forms $\T$,
the stabilizer $\mathrm{G}_{\T}$ is $2$-dimensional, i.e., both parameters  
$k_2, \, k_3 = 0$ cannot vanish. Consequently,  the
holonomy representation $\mathrm{G}_{\T^{\mathrm{c}}} \rightarrow
\U(3)$ splits into the sum of two faithful representations. 
The results of \cite[Lemma 4.4 and Lemma 5.6]{CleytonSwann} yield again that 
the curvature $\mathrm{R}^{\mathrm{c}}$ of the characteristic
connection is $\nabla^{\mathrm{c}}$-parallel,
$\nabla^{\mathrm{c}} \T^{\mathrm{c}} = 
\nabla^{\mathrm{c}} \mathrm{R}^{\mathrm{c}} = 0 \, .$ Since the group
$\mathrm{G}_{\T}$ is compact, the Ambrose-Singer manifold is regular, i.e.,
the manifold  $M^6 = \mathrm{G}/\mathrm{G}_{\T^{\mathrm{c}}}$ is homogeneous and
the Lie algebra of its automorphism group $\mathrm{G}$ 
is the vector space
$\g := \g_{\T^{\mathrm{c}}} \oplus \R^6$ equipped with the bracket
(see \cite[Theorem 5.10]{CleytonSwann})
\bdm
\big[ A \, + \, X \, , \, B \, + \, Y \big] \ = \ 
- \, \mathrm{R}^{\mathrm{c}}(X , Y) \, + \, 
\big( A \cdot Y \, - \, B \cdot X \, - \, \T^{\mathrm{c}}(X , Y) \big) \, . 
\edm
The curvature operator $ \mathrm{R}^{\mathrm{c}} : \Lambda^2(\R^6)
\rightarrow \g_{\T^{\mathrm{c}}}$ is invariant. 
Fix an element
$\omega \in \g_{\T^{\mathrm{c}}}$ of length one and denote by
$\mathrm{R}_{ij}$ the coefficients of the curvature,
$\mathrm{R}^{\mathrm{c}}(e_i \wedge e_j) :=  \mathrm{R}_{ij} \cdot \omega$.
Let us introduce the following element inside the Clifford algebra,
\bdm
\mathrm{R}^{\mathrm{c}} \ := \ \sum_{i<j}\mathrm{R}_{ij} \cdot e_i \cdot e_j 
\cdot \omega \ . 
\edm 
The Jacobi identity implies that
the sum $(\T^{\mathrm{c}})^2 \, + \,\mathrm{R}^{\mathrm{c}}$ 
is a scalar in the Clifford algebra $\mathfrak{Cliff}(\R^6)$ 
(vice versa: if $\mathrm{R} : \Lambda^2(\R^6) \rightarrow \g_{\T^{\mathrm{C}}}
\subset \Lambda^2(\R^6)$ is symmetric, then 
$(\T^{\mathrm{c}})^2 \, + \,\mathrm{R} \in \R^1$ is equivalent
to the Jacobi identity). 
This system
of equations links the curvature operator to the torsion form. We again use
a suitable matrix representation of the Clifford algebra in order to
discuss the system for concrete $3$-forms $\T$.
\begin{lem} \label{LemmaAlgebra} 
There is no $5$-tuple consisting of three
integers $k_1, k_2,  k_3 \in \Z$, a $3$-form
$\T \in \Lambda^3_{12}$ and a curvature operator 
$\mathrm{R}$ such that
\begin{enumerate}
\item $\g_{\T} = \R^1 \cdot (k_1 \cdot e_{12} \, + \, k_2 \cdot e_{34} \, + \, k_3 \cdot e_{56})$ is $1$-dimensional.  
\item The element 
$\T^2 \, + \, \mathrm{R}$ is a scalar   
in $\mathfrak{Cliff}(\R^6)$.
\end{enumerate}
\end{lem}
\begin{proof} 
Since the isotropy algebra $\g_{\T}$ is $1$-dimensional,
the $3$-form is not zero. The $2$-form $(k_1 \cdot e_{12}  +  k_2 \cdot e_{34} 
+  k_3 \cdot e_{56})$ preserves a non-trivial element 
in $\Lambda^3_{12}$ if and only if
\bdm
k_1 \, k_2 \, k_3 \, (k_1 + k_2 - k_3) \, (k_1 - k_2 + k_3) \, 
(- k_1 + k_2 + k_3) \ = \ 0 \, . 
\edm
Basically, there are two cases to consider: that one of the $k_i$'s is zero, $k_3 = 0$,
or that $ k_3 = k_1 + k_2$.\\
Case $1$: $k_3 = 0$. The equation $\rho_*(k_1\cdot e_{12}+ k_2\cdot
e_{34})\T\,=\,0$ reads as
\begin{eqnarray*}
&&k_2 \, A_1 \ = \ k_2 \, A_2 \ = \ k_1 \, A_5 \ = \ k_1 \, A_6 \ 
= \ (k_1 + k_2) \, A_{11} \ = \ (k_1 + k_2) \, A_{12} \ = \ 0 \ ,\\
&&(k_1 - k_2) \, (A_8 \pm A_9) \ = \ 
(k_1 - k_2) \, (A_7 \pm A_{10}) \ = \ 0 \, .
\end{eqnarray*}
We split the first case into four sub-cases:\\
Case $1.1$: $k_1 \neq 0 \neq k_2, \, k_1 + k_2 \neq 0, \, k_1 - k_2 \neq
0, \, k_3 = 0$. The solution space is $2$-dimensional and parameterized 
by the parameters $A_3, \, A_4$ of the $3$-form $\T \in \Lambda^3_{12}$. 
Any $\T$ of that type is preserved by two elements of the Lie algebra $\un(3)$,
$\rho_*(e_{12}) \T = 0 = \rho_*(e_{34}) \T$, hence the dimension of the 
isotropy algebra equals two, a contradiction.\\
Case $1.2$:  $k_1 = k_3 = 0, \, k_2 \neq 0$. The solution space
 $\rho_*(k_1 \cdot e_{12} \, + \, k_2 \cdot e_{34} \, + \, k_3 \cdot e_{56}) 
\T = 0$ is $4$-dimensional and any of these $3$-forms $\T$ 
has a $2$-dimensional
isotropy algebra $\g_{\T}$.\\
Case $1.3$: $k_1 \neq 0 \neq k_2 , \, k_3 = 0 , \, k_1 - k_2 = 0$. 
The solution space is $6$-dimensional and parameterized by the parameters
$A_3, \, A_4 , \, A_7 , \, A_8 , \, A_9 , \, A_{10}$.
For any of these forms, we compute the endomorphism $\T^2 + \mathrm{R}$
in  $\mathfrak{Cliff}(\R^6) \subset \mathfrak{Cliff}(\R^7) \rightarrow
\mathrm{End}(\Delta_7)$ in the $7$-dimensional spin representation. The
condition that $\T^2 + \mathrm{R}$ should be a scalar leads to the following
restrictions
\bdm
\mathrm{R}_{15} \, = \,\mathrm{R}_{16} \, = \, \mathrm{R}_{25} \, = \,  
\mathrm{R}_{26} \, = \, \mathrm{R}_{35} \, = \,  \mathrm{R}_{36} \, = \,
\mathrm{R}_{45} \, = \, \mathrm{R}_{46} \, = \, \mathrm{R}_{56} \, = \, 0 \ ,
\edm
\bdm
\mathrm{R}_{23} \, = \, - \, \mathrm{R}_{14} \, , \quad
\mathrm{R}_{24} \, = \, \mathrm{R}_{13} \ ,
\edm
\bdm
\mathrm{R}_{12} \, = \, - \, \mathrm{R}_{34} \, - \, 2 \sqrt{2} \cdot 
\big( A_3^2 \, + \, A_4^2 \, + \, A_7^2 \, + \, A_8^2 \, + \, A_9^2 \, + \,
A_{10}^2 \big) \ .
\edm
Moreover, the coefficients of the $3$-form have to satisfy the three relations
\bdm
A_3 \, A_{10} \ = \ - \, A_4 \, A_9 \, , \quad
A_7 \, A_{10} \ = \ - \, A_8 \, A_9 \, , \quad
A_3 \, A_8 \ = \ A_4 \, A_7 \ .
\edm
The isotropy algebra $\g_{\T}$ of any $3$-form satisfying these
conditions has dimension two,
i.e., case 1.3 is impossible.\\
Case $1.4$: $ k_1 \neq 0 \neq k_2 , \, k_3 = 0 , \, k_1 + k_2 = 0$. 
The solution space is $4$-dimensional and parameterized 
by  $A_3, \, A_4 , \, A_{11} , \, A_{12}$. The condition
$\T^2 + \mathrm{R} \in \R^1$ for some curvature operator implies 
in particular that two
of the parameters of the $3$-form vanish, 
$A_{11} \, = \, A_{12} \, = \, 0$.
This family of forms has been investigated already in Case $1.1$.
We obtain $\mathrm{dim}(\g_{\T}) = 2$, a contradiction.\\
Case $2$: $k_1 \, k_2 \, k_3 \neq 0 , \, k_3 = k_1 + k_2$. The second
case is simpler. We solve again the equation
$\rho_* \big(k_1 \cdot e_{12} + k_2 \cdot e_{34} + (k_1 + k_2) \cdot e_{56}
\big) \T = 0$. The solution space is $2$-dimensional and parameterized by
the parameters $A_{11}, A_{12}$. Any of these forms has a 
$4$-dimensional isotropy algebra, again a contradiction.
\end{proof}
\noindent
A direct consequence of Lemma \ref{LemmaAlgebra} is the following 
\begin{thm} 
Complete hermitian manifolds $(M^6, g, \J)$ 
of type $\mathcal{W}_3$ such that
\bdm
\nabla^{\mathrm{c}} \T^{\mathrm{c}} \ = \ 0 \, , \quad 
\mathrm{dim}\big(\mathrm{G}_{\T^{\mathrm{c}}}\big) \ = \ 1 
\edm
do not exist.
\end{thm}
\noindent
Consider hermitian  $\mathcal{W}_3$-manifolds 
$(M^6 , g , \J)$ with parallel characteristic torsion and 
$2$-dimensional isotropy group,
\bdm
\nabla^{\mathrm{c}} \T^{\mathrm{c}} \ = \ 0 \, , \quad 
\mathrm{G}_{\T^{\mathrm{c}}} \ 
= \ \mathrm{S}^1 \times \mathrm{S}^1 .
\edm
The curvature of such a hermitian structure is not necessarily parallel,
i.e., $M^6$ does not have to be homogeneous.
Naturally reductive hermitian manifolds can be constructed 
out of a $3$-from $\T \in \Lambda^2_{12}$  and a curvature 
tensor $\mathrm{R} : \Lambda^2(\R^6) \rightarrow \g_{\T}$ 
such that the pair
$(\T \, , \mathrm{R})$ defines a Lie algebra structure on
$\g := \g_{\T} \oplus \R^6$. The naturally reductive space
$\mathrm{G}/\mathrm{G}_{\T}$ is a hermitian $6$-manifold of type 
$\mathcal{W}_3$ with parallel characteristic torsion $\T$. 
\begin{exa} 
The isotropy algebra of the form $\T := e_{125} -  e_{345}$
is generated by $\omega_1 := e_{12}$ and $\omega_2 := e_{34}$. The most
general invariant $2$-form with values in $\g_{\T}$ depends on $6$ parameters,
\bdm
\mathrm{R} \ := \ \sum_{k=1}^2 \big( \mathrm{R}_{12}^k \cdot e_1 \wedge e_2 \,
+ \, \mathrm{R}_{34}^k \cdot e_3 \wedge e_4 \,
+ \, \mathrm{R}_{56}^k \cdot e_5 \wedge e_6 \big) \otimes \omega_k \ . 
\edm
The Jacobi identity is equivalent to 
\bdm
\mathrm{R}_{56}^1 \ = \ \mathrm{R}_{56}^2 \ = \ 0 \, , \quad
\mathrm{R}_{12}^2 \ = \ \mathrm{R}_{34}^1 \ = \ - \, 1 \, .
\edm
There exists a $2$-parameter family of curvature operators associated 
with the form $\T$,
\bdm
\mathrm{R} \ = \ ( \mathrm{R}_{12}^1 \cdot e_1 \wedge e_2 \,
- \, e_3 \wedge e_4 ) \otimes \omega_1 \ + \ 
(- \, e_1 \wedge e_2 \, + \,\mathrm{R}_{34}^2 \cdot e_3 \wedge e_4 )
\otimes \omega_2 \ .
\edm
The holonomy algebra $\h$ of the connection is $1$-dimensional if and only if
$\mathrm{R}_{12}^1 \cdot \mathrm{R}_{34}^2 = 1$ holds. The Lie algebra
$\g = \g_{\T} \oplus \R^6$ has a $2$-dimensional center,
\bdm
\z \ = \ \mathrm{Lin}(\omega_1  -  \omega_2  +  e_5 \, , \, e_6) \, .
\edm
Consider  the Lie algebra $\g^* := \g/\z$. Then $\g$ is a central extension of
$\g^*$. The projections into $\g^*$ of the elements $\omega_1 , \omega_2 , 
e_1 , e_2 , e_3 , e_4$ form a basis of the vector space $\g^*$ and the
commutator relations in $\g^*$ are given by the formulas
\begin{eqnarray*}
{[\omega_1 \, , \, \omega_2 ] }&=& {0 \, , \quad [ \omega_1 \, , \, e_1 ] \ = \ e_2
\, , \quad [ \omega_1 \, , \, e_2 ] \ = \ - \, e_1 \, , 
 \quad [ \omega_1 \, , \, e_3 ] \ = \ [ \omega_1 \, , \, e_4] \ = \ 0 \, , }\\ 
{[ \omega_2 \, , \, e_1 ]} &=& {[\omega_2 \, , \, e_2 ] \ = \ 0 \, , \quad
 [ \omega_2 \, , \, e_3 ] \ = \ e_4 \, , \quad [ \omega_2 \, , \, e_4 ] \ = \
 - \, e_3 \, , } \\
{[ e_1 \, , \, e_3 ] } &=& {[ e_1 \, , \, e_4 ] \ = \ [ e_2 \, , \, e_3 ] 
\ = \ [ e_2 \, , \, e_4 ] \ = \ 0 \, , }\\
{ [ e_1 \, , \, e_2 ] } &=& {(1 - \mathrm{R}_{12}^1) \, \omega_1 \, , \quad 
[ e_3 \, , \, e_4 ] \ = \ (1 - \mathrm{R}_{34}^2) \, \omega_2  } \ .
\end{eqnarray*}
The Lie algebra $\g^*$ is the sum of 
two subalgebras,
\bdm
\p_1 \ = \ \Lin \big( \omega_1 \, , \, e_1 \, , \, e_2 \big) \, , \quad
\p_2 \ = \ \Lin \big( \omega_2 \, , \, e_3 \, , \, e_4 \big) \, , 
\edm
and we have
\bdm
\g^* \ = \ \p_1 \oplus \p_2 \, , \quad
[\p_1 \, , \, \p_1 ] \ \subset \p_1 \, , \quad
[\p_2 \, , \, \p_2 ] \ \subset \p_2 \, , \quad
[\p_1 \, , \, \p_2 ] \ = \ 0. 
\edm
The Lie algebras $\p_1 \, , \, \p_2$ are isomorphic to $\so(3,\R)$,
$\slin(2,\R)$ or to the $3$-dimensional nilpotent Lie algebra.
Consequently, we gave a complete description of the possible
automorphism groups of all naturally reductive hermitian
$\mathcal{W}_3$-manifolds with parallel characteristic
torsion of type $\T = e_{125} - e_{345}$. 
\end{exa}
\noindent
We remark that the torsion form $\T = e_{125} - e_{345}$ represents the
general case. Indeed, let $\T \in \Lambda_{12}^3$ be a $3$-form with
a $2$-dimensional isotropy group. Following once again carefully
the proof of lemma \ref{LemmaAlgebra}, we see that this form behaves like
$e_{125} - e_{345}$ in the sense that the automorphism groups are
the same. Therefore, we obtain:
\begin{thm}
Any naturally reductive hermitian $\mathcal{W}_3$-manifold with
a $2$-dimensional isotropy algebra $\g_{\T^{\mathrm{c}}}$ of its characteristic
torsion is locally isometric to one of the spaces described in the
previous example.
\end{thm}
%
%
    

\begin{thebibliography}{1111}
\bibitem[1]{AGS}
E. Abbena, S. Garbiero and S. Salamon, \emph{Almost hermitian
geometry on six dimensional nilmanifolds}, Ann. Sc. Norm. Sup. 30 (2001),
147-170.
\bibitem[2]{Adams}
J.F. Adams, \emph{Lectures on Lie groups}, University of Chicago Press 1969.
\bibitem[3]{Agri}
I. Agricola, \emph{Connections on naturally reductive spaces, their
Dirac operator and homogeneous models in string theory}, Comm. 
Math. Phys. 232 (2003), 535-563.
\bibitem[4]{AgFr1}
I. Agricola and Th. Friedrich, \emph{On the holonomy of connections
with skew-symmetric torsion}, Math. Ann. 328 (2004), 711-748.
\bibitem[5]{AgFr2}
I. Agricola and Th. Friedrich, \emph{The Casimir operator of a metric 
connection with skew-symmetric torsion}, Journ. Geom. Phys. 50 (2004), 188-204.
\bibitem[6]{Alexandrov}
B. Alexandrov, \emph{$\Sp(n)\U(1)$-connections with parallel totally 
skew-symmetric torsion}, math.dg/0311248.
\bibitem[7]{DAtriZiller}
J.E. D'Atri and W. Ziller, \emph{Naturally reductive metrics and Einstein
metrics on compact Lie groups}, Memoirs of AMS no. 215, 1979.
\bibitem[8]{BFGK}
H. Baum, Th. Friedrich, R. Grunewald and I. Kath, \emph{Twistors and
Killing spinors on Riemannian manifolds}, Teubner-Texte zur Mathematik
No. 124, Teubner-Verlag Leipzig/Stuttgart 1991.
\bibitem[9]{Belgun}
F. Belgun, \emph{On the metric structure of non-K\"ahler complex surfaces}, 
Math. Ann. 317 (2000), 1-40.
\bibitem[10]{BelgunMoroianu}
F. Belgun and A. Moroianu, \emph{Nearly K\"ahler $6$-manifolds with
reduced holonomy}, Ann. Global Anal. Geom. 19 (2001), 307-319.
\bibitem[11]{Bismut}
J. M. Bismut, \emph{A local index theorem for non-K\"ahlerian manifolds},
Math. Ann. 284 (1989), 681-699.
\bibitem[12]{Butruille}
J.-B. Butruille, \emph{Classification des varietes approximativement
k\"ahleriennes homogenes}, math.dg/0401152.
\bibitem[13]{ChiossiSalamon}
S. Chiossi and S. Salamon, \emph{The intrinsic torsion of $\SU(3)$ and
$\G_2$-structures}, Differential geometry, Valencia 2001, 115-133. Word Sci.
Publishing, River Edge, NJ., 2002.
\bibitem[14]{CleytonSwann}
R. Cleyton and A. Swann, \emph{Einstein metrics via intrinsic or 
parallel torsion}, math.dg/0211446.
\bibitem[15]{FaFaSalamon}
M. Falcitelli, A. Farinola and S. Salamon, \emph{Almost-hermitian geometry},
Differential Geom. Appl. 4 (1994), 259-282.
\bibitem[16]{FinoPortonSalamon}
A. Fino, M. Porton and S. Salamon, \emph{Families of strong KT structures in
six dimensions}, math.dg/0209259.
\bibitem[17]{Fri1}
Th. Friedrich, \emph{Der erste Eigenwert des Dirac Operators einer 
kompakten Riemannschen Mannigfaltigkeit nichtnegativer 
Skalarkr\"ummung},
Math. Nachr. 97 (1980), 117-146. 
\bibitem[18]{Fri2}
Th. Friedrich, \emph{Dirac Operators in Riemannian Geometry}, 
Graduate Studies in Mathematics vol. 25, AMS, Providence 2000.
\bibitem[19]{Fri3}
Th. Friedrich, \emph{On types of non integrable geometries}, 
Rend. Circ. Mat. di Palermo 71 (2003), 99-113.
\bibitem[20]{Fri4}
Th. Friedrich, \emph{$\mathrm{Spin}(9)$-structures and connections
with totally skew-symmetric torsion}, Journ. Geom. Phys. 47 (2003), 197-206.
\bibitem[21]{FriedrichGrunewald}
Th. Friedrich and R. Grunewald, \emph{On the first eigenvalue
of the Dirac operator on $6$-dimensional manifolds}, Ann. Global
Anal. Geom. 3 (1985), 265-273.
\bibitem[22]{FriedrichIvanov}
Th. Friedrich and S. Ivanov, \emph{Parallel spinors and connections with
skew-symmetric torsion in string theory}, Asian Journ. Math. 6 (2002), 
303-336.
\bibitem[23]{FriedrichKurke}
Th. Friedrich and H. Kurke, \emph{Compact four-dimensional self-dual Einstein
manifolds with positive scalar curvature}, Math. Nachr. 106 (1982), 
271-299.
\bibitem[24]{GrantcharovPonn}
D. Grantcharov, G. Grantcharov and Y.S. Poon, \emph{Calabi-Yau connections
with torsion on toric bundles}, math.dg/0306207.
\bibitem[25]{Gray69}
A. Gray, \emph{Almost complex submanifolds of the six sphere}, 
Proc. Amer. Math. Soc. 20 (1969), 277-279.
\bibitem[26]{Gray69a}
A. Gray, \emph{Six-dimensional almost complex manifolds defined by
means of the three-fold vector cross products}, Tohoku Math. Journ. 
II. Ser., 21 (1969), 614-620. 
\bibitem[27]{Gray70}
A. Gray, \emph{Nearly K\"ahler manifolds}, Journ. Diff. Geom. 4 (1970),
283-310. 
\bibitem[28]{Gray76}
A. Gray, \emph{The structure of nearly K\"ahler manifolds}, Math. Ann. 223 (1976), 233-248.
\bibitem[29]{GrayHervella}
A. Gray and L. Hervella, \emph{The sixteen classes of almost
hermitian manifolds and their linear invariants}, Ann. di Mat.
Pura ed Appl. 123 (1980), 35-58.
\bibitem[30]{Grunewald}
R. Grunewald, \emph{Six-dimensional Riemannian manifolds with
real Killing spinors}, Ann. Glob. Anal. Geom. 8 (1990), 43-59.
\bibitem[31]{Hitchin74}
N. Hitchin, \emph{On compact four-dimensional Einstein manifolds}, 
Journ. Diff. Geom. 9 (1974), 435-442.
\bibitem[32]{Hitchin}
N. Hitchin, \emph{K\"ahlerian twistor spaces}, Proc. Lond. Math. Soc.
III Ser., 43 (1981), 133-150.
\bibitem[33]{Kirichenko}
V. Kirichenko, \emph{$\mathrm{K}$-spaces of maximal rank}, Mat. Zam. 
22 (1977), 465-476.
\bibitem[34]{Kostant}
B. Kostant, \emph{A cubic Dirac operator and the emergence of Euler number
multiplets of representations for equal rank subgroups}, Duke Math. J. 100
(1999), 447-501. 
\bibitem[35]{Matsumoto}
M. Matsumoto, \emph{On $6$-dimensional almost Tachibana spaces}, 
Tensor N. S. 23 (1972), 250-282.
\bibitem[36]{Nagy}
P.-A. Nagy, \emph{Nearly K\"ahler geometry and Riemannian foliations},
Asian J. Math. 6 (2002), 481-504.
\bibitem[37]{Sternberg}
S. Sternberg, \emph{Lie algebras}, preprint November 1999. 
\bibitem[38]{Takamatsu}
K. Takamatsu, \emph{Some properties of $6$-dimensional $\mathrm{K}$-spaces}, 
Kodai Math. Sem. Rep. 23 (1971), 215-232.
\bibitem[39]{Vaisman}
I. Vaisman, \emph{Locally conformal K\"ahler manifolds with parallel Lee
  form}, Rendiconti di Matem., Roma, 12 (1979), 263-284.
\bibitem[40]{YanoKon}
K. Yano and M. Kon, \emph{Structures on manifolds}, World Scientific 1984. 
\end{thebibliography}
\end{document}